\documentclass[11pt]{amsart}

\usepackage{amsmath,amssymb,amsthm,amsfonts,amscd,flafter,epsf, epsfig,graphicx,verbatim,pinlabel,mathrsfs}
\usepackage[all]{xy}
\usepackage{epsf} 

\title{Capping off open books and the Ozsv{\'a}th-Szab{\'o} contact invariant}
\author[John A. Baldwin]{John A. Baldwin}
\address {Department of Mathematics, Princeton University\\ Princeton, NJ 08544-1000}
\email {baldwinj@math.princeton.edu}
\thanks{The author was partially supported by an NSF Postdoctoral Fellowship.}
\date{}

\newcommand\bx{\mathbf{x}}
\newcommand\bF{\mathbb{F}}
\newcommand\bR{\mathbf{R}}
\newcommand\by{\mathbf{y}}

\newcommand\bn{\mathbf{n}}
\newcommand\bw{\mathbf{w}}

\newcommand\bt{\mathbf{\Theta}}

\newcommand\hf{\widehat{HF}}
\newcommand\hfp{HF^+}

\newcommand\cfp{CF^+}

\newcommand\zz{\mathbb{Z}}
\newcommand\zzt{\mathbb{Z}_2}

\newcommand\Sc{\text{Spin}^c}
\newcommand\spc{\mathfrak{s}}
\newcommand\spt{\mathfrak{t}}
\newcommand\Ta{\mathbb{T}_{\alpha}}
\newcommand\Tb{\mathbb{T}_{\beta}}
\newcommand\Tc{\mathbb{T}_{\gamma}}

\newtheorem{theorem}{Theorem}[section]
\newtheorem{lemma}[theorem]{Lemma}

\newtheorem{conjecture}[theorem]{Conjecture}
\newtheorem{corollary}[theorem]{Corollary}
\newtheorem{proposition}[theorem]{Proposition}
\newtheorem{question}[theorem]{Question}

\theoremstyle{definition}

\newtheorem{remark}[theorem]{Remark}

\newtheorem{example}[theorem]{Example}

\begin{document}
\begin{abstract} 
If $(S,\phi)$ is an open book with disconnected binding, then we can form a new open book $(S',\phi')$ by capping off one of the boundary components of $S$ with a disk. Let $M_{S,\phi}$ denote the 3-manifold with open book decomposition $(S,\phi)$. We show that there is a $U$-equivariant map from $\hfp(-M_{S',\phi'})$ to $\hfp(-M_{S,\phi})$ which sends $c^+(S',\phi')$ to $c^+(S,\phi)$, and we discuss various applications. In particular, we determine the support genera of almost all contact structures which are compatible with genus one, one boundary component open books. In addition, we compute $d_3(\xi)$ for every tight contact manifold $(M,\xi)$ supported by a genus one open book with periodic monodromy.
\end{abstract}

\maketitle

\section{Introduction}
Giroux's correspondence between contact structures up to isotopy and open books up to positive stabilization allows us to translate questions from contact geometry into questions about diffeomorphisms of compact surfaces with boundary \cite{giroux2}. As a result, one is inclined to wonder about the contact-geometric significance of certain natural operations which can be performed on open books. For instance, let us define the ``composition" of two open books $(S,\phi_1)$ and $(S,\phi_2)$ to be the open book $(S,\phi_1 \circ \phi_2)$. In \cite{bald3}, we use the Ozsv{\'a}th-Szab{\'o} contact invariant (see \cite{osz1}) to study the relationship between the contact structures supported by two such open books and the contact structure supported by their composition. There, we prove the following.
\begin{theorem}[{\rm \cite{bald3}}]
\label{thm:comult}
If $c(S,\phi_1)$ and $c(S,\phi_2)$ are both non-zero, then so is $c(S,\phi_1 \circ \phi_2)$.
\end{theorem}
In particular, Theorem \ref{thm:comult} implies that if the contact structures supported by $(S,\phi_1)$ and $(S,\phi_2)$ are both strongly fillable, then the contact structure supported by $(S,\phi_1 \circ \phi_2)$ is tight (forthcoming work by the author \cite{bald8} and, independently, by Baker, Etnyre and Van Horn-Morris \cite{bev2} strengthens this result).

In this paper, we use the Ozsv{\'a}th-Szab{\'o} contact invariant to study the geometric effect of another natural operation on open books called ``capping off". Consider the open book $(S_{g,r},\phi)$, where $S_{g,r}$ is a genus $g$ surface with $r>1$ boundary components and $\phi$ is a diffeomorphism of $S_{g,r}$ which fixes the boundary pointwise. By capping off one of the boundary components of $S_{g,r}$ with a disk, we obtain an open book $(S_{g,r-1},\phi'),$ where $\phi'$ is the extension of $\phi$ to $S_{g,r-1}$ by the identity on this disk. Let $M_{S,\phi}$ denote the 3-manifold with open book decomposition $(S,\phi)$. There is a natural cobordism $W$ from $M_{S_{g,r},\phi}$ to $M_{S_{g,r-1},\phi'}$ obtained by attaching a $0$-framed 2-handle along the binding component in $M_{S_{g,r},\phi}$ corresponding to the capped off boundary component of $S_{g,r}$. Alternatively, we can think of $W$ as a cobordism from $-M_{S_{g,r-1},\phi'}$ to $-M_{S_{g,r},\phi}$. Our main theorem is the following.

\begin{theorem}
\label{thm:nat}
There exists a $\Sc$ structure $\spc_0$ on $W$ for which the map $$F^+_{W,\spc_0}:\hfp(-M_{S_{g,r-1},\phi'})\rightarrow \hfp(-M_{S_{g,r},\phi})$$ sends $c^+(S_{g,r-1},\phi')$ to $c^+(S_{g,r},\phi).$
\end{theorem}

In some sense, Theorem \ref{thm:nat} generalizes Ozsv{\'a}th and Szab{\'o}'s original definition of the contact invariant. Their definition begins with the fact that every contact structure is supported by an open book of the form $(S_{g,1},\phi)$ for some $g>1$. By capping off the one boundary component of $S_{g,1}$, we obtain a closed surface $S_{g,0}$. If $M_{S_{g,0},\phi'}$ is the corresponding fibered 3-manifold and $\mathfrak{t}$ is the $\Sc$ structure on $M_{S_{g,0},\phi'}$ represented by the vector field transverse to the fibers, then $\hfp(-M_{S_{g,0},\phi'}, \mathfrak{t})$ is generated by a single element $c^+(S_{g,0},\phi')$ \cite{osz17}. In \cite{osz1}, Ozsv{\'a}th and Szab{\'o} define the contact invariant $c^+(S_{g,1},\phi)$ to be the image of this element $c^+(S_{g,0},\phi')$ under the map $$F^+_V:\hfp(-M_{S_{g,0},\phi'}) \rightarrow \hfp(-M_{S_{g,1},\phi})$$ induced by the corresponding 2-handle cobordism $V$.

At first glance, this 2-handle attachment does not seem like a very natural contact-geometric operation, though Eliashberg proves in \cite{yasha2} that there is a symplectic form $\Omega$ on the cobordism $V$ which is positive on the fibers of the fibration $M_{S_{g,0},\phi'}\rightarrow S^1$, and for which the contact 3-manifold supported by $(S_{g,1},\phi)$ is a weakly concave boundary component of $(V,\Omega)$. One expects that a similar construction should produce a symplectic structure on the cobordism $W$ considered in Theorem \ref{thm:nat}. In fact, since this paper first appeared, Gay and Stipsicz have shown that one can find a symplectic form on $W$ for which the contact 3-manifolds supported by $(S_{g,r},\phi)$ and $(S_{g,r-1},\phi')$ are strongly concave and strongly convex, respectively \cite{gays}. On the other hand, the contact invariant in Heegaard Floer homology (in contrast with its analogue in Monopole Floer homology \cite{mr}) is not known, in general, to behave naturally with respect to the map induced by a strong symplectic cobordism, and so Theorem \ref{thm:nat} provides new information in this regard.

Below, we explore some consequences and potential applications of Theorem \ref{thm:nat}, and we discuss some natural questions which arise from this result. To begin with, consider the following immediate corollary of Theorem \ref{thm:nat}.

\begin{corollary}
\label{cor:OT}
If $c^+(S_{g,r-1},\phi') = 0$ then $c^+(S_{g,r},\phi) = 0$.
\end{corollary}

This prompts the question below.

\begin{question}
\label{question:OT}
Is the contact structure supported by $(S_{g,r},\phi)$ is overtwisted whenever the contact structure supported by $(S_{g,r-1},\phi')$ is?
\end{question}

A positive answer to this question would be helpful in determining which open books can support tight contact structures. For example, the genus one, one boundary component open books which support tight contact structures are classified in \cite{bald1,hkm3, hkm2}. Combined with this classification, Corollary \ref{cor:OT} (or, an affirmative answer to Question \ref{question:OT}) places concrete restrictions on genus one open books with multiple boundary components which can support tight contact structures.

If $W$ is the 2-handle cobordism in Theorem \ref{thm:nat}, then the map $F^+_W$ (summing over all $\Sc$ structures) fits into a surgery exact triangle, 
$$\xymatrix @!0@C=5pc@R=5pc{
  \hfp(-M_{S_{g,r-1},\phi'}) \ar[rr]^-{F^+_W}&  &\hfp(-M_{S_{g,r},\phi}) \ar[dl]^<(.36){F^+_{X}} \\
 &\hfp(-M_{S_{g,r},\phi\cdot t_{\gamma}^{-1}}), \ar[ul] }$$ where $t_{\gamma}^{-1}$ is a left-handed Dehn twist around a curve $\gamma$ parallel to the boundary component $B$ that we are capping off \cite{osz14}. In this triangle, $X$ is the cobordism obtained by attaching a $(-1)$-framed 2-handle to $-M_{S_{g,r},\phi}$ along the binding component corresponding to $B$. According to the following theorem, the contact invariant behaves naturally under the map induced by $X$ as well (strictly speaking, the theorem below is only proved in \cite{hkm3} for curves $\gamma$ which are non-separating, but as long as $r>1$, we can stabilize the open book and then apply the result from \cite{hkm3}).
 
 \begin{theorem}[{\rm \cite{osz1,hkm3}}]
 \label{thm:nat2}
The map $F^+_X$ sends $c^+(S_{g,r},\phi)$ to $c^+(S_{g,r},\phi \cdot t_{\gamma}^{-1}).$
\end{theorem}

The surgery exact triangle has proved to be one of the most versatile tools in Heegaard Floer homology. One therefore expects that the triangle above, combined with Theorems \ref{thm:nat} and \ref{thm:nat2}, may be used to provide interesting contact-geometric information in many settings. For example, we obtain the following obstruction to there being a Stein structure on $W$ (it is clear that the 2-handle attachment used to form $W$ cannot locally be done in a Stein way; however, it may sometimes be possible to globally construct such a Stein structure).

\begin{corollary}
\label{cor:stein}
The cobordism $W:M_{S_{g,r},\phi} \rightarrow M_{S_{g,r-1},\phi'}$ has a Stein structure compatible with the contact structures on either end only if $c^+(S_{g,r},\phi\cdot t_{\gamma}^{-1}) = 0$.
\end{corollary}

For, if $W$ has a Stein structure, then the map $F^+_W$ sends $c^+(S_{g,r-1},\phi')$ to $c^+(S_{g,r},\phi)$ (see \cite{osz1}), and the exactness of the triangle above implies that $$c^+(S_{g,r},\phi\cdot t_{\gamma}^{-1}) = F^+_X(c^+(S_{g,r},\phi)) = 0.$$

Below, we describe a consequence of Theorem \ref{thm:nat} for contact surgery on stabilized Legendrian knots. Suppose $K$ is an oriented Legendrian knot in $(M,\xi)$, and let $(M_{\pm 1}(K),\xi_{\pm 1}(K))$ be the contact 3-manifold obtained from $(M,\xi)$ via contact $\pm 1$-surgery on $K.$ We denote by $S_+(K)$ and $S_-(K)$ the positive and negative Legendrian stabilizations of $K$, as defined in \cite{EH2}. Let $K'$ be either $S_+(K)$ or $S_-(K)$. As we shall see in Section \ref{sec:legstab}, the following is a special case of Theorem \ref{thm:nat}. 

\begin{theorem}
\label{thm:stab}
There is a $U$-equivariant map $F^+:\hfp(-M_{\pm1}(K))\rightarrow \hfp(-M_{\pm 1}(K'))$ which sends $c^+(\xi_{\pm 1}(K))$ to $c^+(\xi_{\pm 1}(K'))$.
\end{theorem}

The operation of capping off is closely related to the operation of gluing two open books together along some proper subset of their binding components. More precisely, suppose that $(S,\phi)$ and $(S',\phi')$ are two open books such that either $S$ or $S'$ has more than $n$ boundary components. Let $B_1,\dots,B_n$ and $B'_1,\dots, B'_n$ denote boundary components of $S$ and $S'$, respectively. One forms a surface $S''$ by gluing $S$ to $S'$ by a map which identifies $B_i$ with $B'_i$ for each $i=1,\dots,n$. And one can define a diffeomorphism $\phi''$ of $S''$ whose restriction to $S\subset S''$ is $\phi$ and whose restriction to $S'\subset S''$ is $\phi''$. We say that $(S'',\phi'')$ is an open book obtained by gluing $(S,\phi)$ to $(S',\phi')$.

\begin{remark}
When $n=1$, the 3-manifold $M_{S'',\phi''}$ is homeomorphic to that corresponding to the \emph{contact fiber sum} of the open books $(S,\phi)$ and $(S',\phi')$  (see \cite{wen} for a recent application of contact fiber sum). In contrast, the contact structure supported by the glued open book $(S'',\phi'')$ is generally different from that associated to the contact fiber sum. \end{remark}

We discuss the relationship between capping off and gluing in more detail in Section \ref{sec:glue}; in particular, we prove the following consequence of Theorem \ref{thm:nat}.

\begin{theorem}
\label{thm:glue}
Suppose that $(S'',\phi'')$ is an open book obtained by gluing $(S,\phi)$ to $(S',\phi')$. If $c(S,\phi)$ and $c(S',\phi')$ are both non-zero, then so is $c(S'',\phi'')$.
\end{theorem}

Our study of the effect of capping off on the Ozsv{\'a}th-Szab{\'o} contact invariant began as an attempt to better understand the \emph{support genus} of contact structures, as defined by Etnyre and Ozbagci in \cite{et3}. The support genus of a contact structure $\xi$ is defined to be the minimum, over all open books $(S,\phi)$ compatible with $\xi$, of the genus of $S$; we denote this invariant by $sg(\xi)$. In 2004, Etnyre showed that all overtwisted contact structures have support genus zero, while there are fillable contact structures with $sg(\xi)>0$ \cite{et}. More recently, Ozsv{\'a}th, Szab{\'o} and Stipsicz have found a Heegaard Floer homology obstruction to $sg(\xi)=0$. Their main result is the following.

\begin{theorem}[{\rm \cite{osz13}}]
\label{thm:planar}
Suppose that $\xi$ is a contact structure on the 3-manifold $M$. If $sg(\xi)=0$, then $c^+(\xi) \in U^d \cdot \hfp(-M)$ for all $d \in \mathbb{N}$.
\end{theorem}

Note that Theorem \ref{thm:planar} follows as an immediate corollary of our Theorem \ref{thm:nat}. For, if $sg(\xi)=0$, then $\xi$ is supported by an open book of the form $(S_{0,r},\phi)$, and we may cap off all but one of the boundary components of $S_{0,r}$ to obtain an open book $(S_{0,1},\phi')$. Since all diffeomorphisms of the disk $S_{0,1}$ are isotopic to the identity, $(S_{0,1},\phi')$ is an open book for $S^3$. And, because every element of $\hfp(S^3)$ is in the image of $U^d$ for all $d \in \mathbb{N}$, Theorem \ref{thm:nat} implies that the same is true of the contact invariant $c^+(S_{0,r},\phi) \in \hfp(-M)$ (since the maps $F^+_{W,\spc_0}$ are $U$-equivariant).

In order to use Theorem \ref{thm:planar} to prove that $sg(\xi)>0$, one must be able to show that $c^+(\xi)$ is not in the image of $U^d$ for some $d \in \mathbb{N}$. In practice, this can be very difficult, though the proposition below is sometimes helpful in this regard. 

\begin{proposition}[{\rm \cite{osz17, osz13}}]
\label{prop:nontor} Suppose that $\xi$ is a contact structure on $M$, and let $\spt_{\xi}$ be the $\Sc$ structure associated to $\xi$. If $c^+(\xi) \neq 0$ and the first Chern class $c_1(\spt_{\xi})$ is non-torsion, then there is some $d \in \mathbb{N}$ for which $c^+(\xi) \notin U^d \cdot \hfp(-M)$; hence $sg(\xi)>0$.
\end{proposition}

Theorem \ref{thm:nat} may be used to extend the reach of this proposition a bit further. This is illustrated by the following example.

\begin{example}
\label{ex:sg} Let $a$, $b$, $c$ and $\gamma$ be the curves on $S_{1,2}$ shown on the left in Figure \ref{fig:example}. If $(M,\xi)$ is the contact 3-manifold supported by the open book $(S_{1,2}, (t_a t_b)^5 t_{\gamma}^2 t_c^2)$, then $c_1(\spt_{\xi})$ is twice a generator of $H_1(M;\zz) \cong \zz$ (see the proof of \cite[Theorem 6.2]{et3}). Moreover, $\xi$ is Stein fillable since the diffeomorphism $(t_a t_b)^5 t_{\gamma}^2 t_c^2$ is a product of right-handed Dehn twists; hence, $c^+(\xi) \neq 0$. By Proposition \ref{prop:nontor}, there is some $d \in \mathbb{N}$ for which $c^+(\xi) \notin U^d \cdot \hfp(-M)$. 

\begin{figure}[!htbp]

\labellist 
\hair 2pt 
\small\pinlabel $B$ at 31 158
\pinlabel $c$ at -9 127
\pinlabel $\gamma$ at 36 50
\pinlabel $b$ at 65 84
\pinlabel $a$ at 102 12

\pinlabel $b$ at 289 84
\pinlabel $a$ at 326 12

\endlabellist 

\begin{center}
\includegraphics[width=8.8cm]{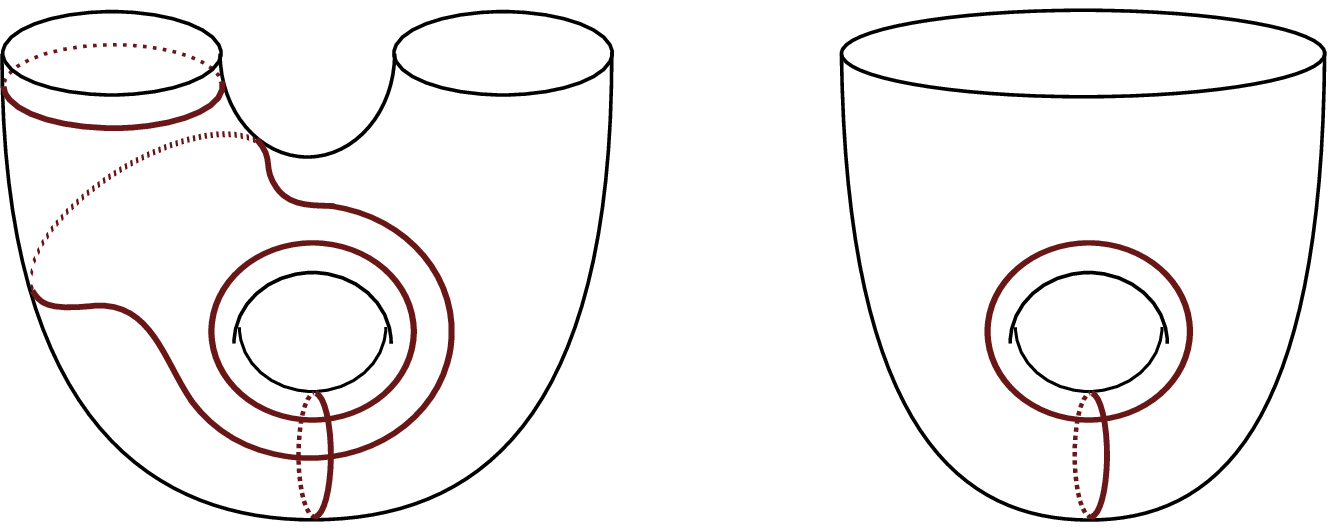}
\caption{\quad The surfaces $S_{1,2}$ and $S_{1,1}$, and the curves $a$, $b$, $c$ and $\gamma$.}
\label{fig:example}
\end{center}
\end{figure}

After capping off the boundary component of $S_{1,2}$ labeled $B$, the curve $c$ becomes null-homotopic, and $\gamma$ becomes isotopic to the curve $b$ on $S_{1,1}$. The capped off open book is therefore $(S_{1,1}, (t_a t_b)^5 t_b^2)$. If $(M',\xi')$ is the contact 3-manifold compatible with this open book, then Theorem \ref{thm:nat} implies that there is some $d \in \mathbb{N}$ for which $c^+(\xi') \notin U^d \cdot \hfp(-M')$ since the same is true for $c^+(\xi)$; hence, $sg(\xi') > 0$. Note that we could \emph{not} have drawn this conclusion directly from Proposition \ref{prop:nontor} since $c_1(\spt_{\xi'})=0$; indeed, the $Spin^c$ structure associated to any contact structure compatible with a genus one, one boundary component open book has trivial first Chern class \cite{et3}. 

\end{example}

In Section \ref{sec:sg}, we use Example \ref{ex:sg} to determine the support genera of almost all contact structures compatible with genus one, one boundary component open books whose monodromies are pseudo-Anosov.

The support genus is not well understood in general, revealing a fundamental gap in our understanding of the link between open books and contact geometry. To begin with, it is not known whether there exist contact structures with $sg(\xi)>1$. Moreover, all of the contact structures that we know of with $sg(\xi)>0$ are at least weakly fillable. It is our hope that Theorem \ref{thm:nat} may be helpful in addressing the first problem. Suppose we wished to find an obstruction to $sg(\xi)=1$. Every genus one open book can be reduced, via capping off, to a genus one open book with one binding component, and much is known about the contact structures compatible with (and the contact invariants associated to) the latter sort of open book \cite{bald1,hkm3, hkm2}. Any property shared by the Ozsv{\'a}th-Szab{\'o} invariants for such contact structures, which is preserved by the map induced by capping off, will provide an obstruction to $sg(\xi)=1$ (we used this principle above to re-derive the obstruction in Theorem \ref{thm:planar} to $sg(\xi)=0$). 

So far, this approach has borne a very modest amount of fruit. A diffeomorphism $\phi$ of $S_{g,r}$ is called \emph{reducible} if $\phi$ is freely isotopic to a diffeomorphism which fixes an essential multi-curve on $S_{g,r}$ (a \emph{free} isotopy is not required to fix points on $\partial S_{g,r}$). We say that $\phi$ is \emph{periodic} if $\phi^m$ is freely isotopic to the identity for some $m \in \mathbb{N}$ and $\phi$ is not reducible. Using the strategy outlined above, we prove the following theorem in Section \ref{sec:capper}.

\begin{theorem}
\label{thm:periodic}
Suppose that the contact 3-manifold $(M,\xi)$ is supported by a genus one open book with $r$ binding components and periodic monodromy. If $\xi$ is tight, then $r\geq -1-4d_3(\xi)$. 
\end{theorem}

Here, $d_3(\xi)$ is the ``3-dimensional" invariant associated to $\xi$, which is well-defined in $\mathbb{Q}$ as long as $c_1(\spt_{\xi})$ is a torsion class. We strengthen Theorem \ref{thm:periodic} at the end of Section \ref{sec:capper}, giving an explicit formula for $d_3(\xi)$ whenever $\xi$ is supported by a genus one open book with periodic monodromy. 

Related to the notion of support genus (and equally mysterious) is that of \emph{binding number} \cite{et3}. If $sg(\xi) = g$, then the binding number of $\xi$ is defined to be the minimum, over all open books $(S_{g,r},\phi)$ compatible with $\xi$, of the number of binding components of the open book, $r$; we denote this invariant by $bn(\xi)$. If $sg(\xi)>0$ and $\xi$ is supported by a genus one open book with periodic monodromy, then Theorem \ref{thm:periodic} implies that $bn(\xi) \geq -1-4d_3(\xi)$. Note that this inequality is sharp for the tight contact structure $\xi_{std}$ on $S^3$, as $bn(\xi_{std}) = 1$ and $d_3(\xi_{std}) = -1/2$.

If $\phi$ is neither reducible nor periodic, then $\phi$ is called \emph{pseudo-Anosov}; these are the most abundant sort. In Section \ref{sec:cappa}, we give a more intrinsic definition of pseudo-Anosov diffeomorphisms, and we discuss properties of these maps which are preserved under capping off. 

\subsection*{Acknowledgements} I wish to thank John Etnyre, Peter Ozsv{\'a}th and Andr{\'a}s Stipsicz for helpful discussions and correspondence.

\section{Proof of Theorem \ref{thm:nat}}
\label{sec:nat}
\subsection{Heegaard diagrams and the contact class}
\label{ssec:heeg}

Let $S$ be a compact surface with boundary, and suppose that $\phi$ is a diffeomorphism of $S$ which restricts to the identity on $\partial S$. Recall that the open book $(S,\phi)$ specifies a closed, oriented 3-manifold $M_{S,\phi}=S \times [0,1] /\sim$, where $\sim$ is the identification given by 
\begin{eqnarray*}
(x,1) \sim (\phi(x),0), && x \in S\\
(x,t) \sim (x, s),  && x\in \partial S, \ t,s \in [0,1].
\end{eqnarray*} 
$M_{S,\phi}$ has a Heegaard splitting $M_{S,\phi} = H_1 \cup H_2$, where $H_1$ is the handlebody $S \times [0,1/2]$ and $H_2$ is the handlebody $S \times [1/2,1]$. Let $S_t$ denote the page $S \times \{t\}.$ The Heegaard surface in this splitting is $$\Sigma := \partial H_1 = S_{1/2} \cup -S_0.$$ If $S=S_{g,r}$ then $\Sigma$ has genus $n=2g+r-1$. To give a pointed Heegaard diagram for $M_{S,\phi}$, it remains to describe the $\alpha$ and $\beta$ attaching curves and the placement of a basepoint $z$.

\begin{figure}[!htbp]
\labellist 
\hair 2pt 
\small\pinlabel $B$ at 50 210
\pinlabel $c_1$ at 5 166 
\pinlabel $a_1$ at 109 112 
\pinlabel $a_{r-1}$ at 367 112 

\pinlabel $a_{r}$ at 456 56 
\pinlabel $a_n$ at 561 56 
\pinlabel $S$ at 260 30 
\endlabellist 
\begin{center}
\includegraphics[width=11.5cm]{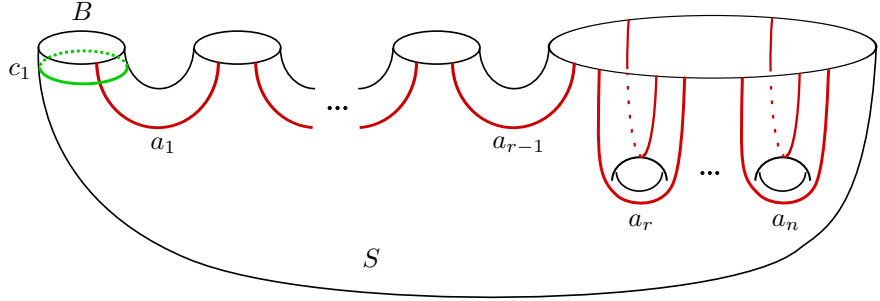}
\caption{\quad The surface $S = S_{g,r}$ and the curve $c_1$. The arcs $a_1, \dots, a_{n}$ are drawn in red.}
\label{fig:page}
\end{center}
\end{figure}

Let $a_1, \dots, a_n$ be pairwise disjoint, properly embedded arcs in $S$ for which the complement $S \setminus\cup a_i$ is a disk. For each $i=1,\dots,n$, let $b_i$ be an arc obtained by changing $a_i$ via a small isotopy which moves the endpoints of $a_i$ along $\partial S$ in the direction specified by the orientation of $\partial S$, so that $a_i$ intersects $b_i$ transversely in one point and with positive sign (where $b_i$ inherits its orientation from $a_i$). For $i=1,\dots,n$, let $\alpha_i$ and $\beta_i$ be the curves on $\Sigma$ defined by $$\alpha_i = a_i \times \{1/2\} \cup a_i \times \{0\},$$ $$\beta_i = b_i \times \{1/2\} \cup \phi(b_i) \times \{0\}.$$ Place a basepoint $z$ in the ``big" region of $S_{1/2}\setminus\cup \alpha_i \setminus\cup \beta_i$ (that is, outside of the thin strip regions), and let $\alpha = \{\alpha_1,\dots,\alpha_n\}$ and $\beta = \{\beta_1,\dots,\beta_n\}$. We say that $(\Sigma, \alpha, \beta, z)$ is a \emph{standard} pointed Heegaard diagram for the open book $(S,\phi)$. See Figure \ref{fig:page} for an example.

\begin{figure}[!htbp]

\labellist 
\hair 2pt 
\small\pinlabel $S$ at 205 35
\pinlabel $S_{1/2}$ at 830 35
\pinlabel $-S_0$ at 830 490
\tiny\pinlabel $\bullet z$ at 710 40 
\pinlabel $x$ at 194 133
\endlabellist 

\begin{center}
\includegraphics[width=11.5cm]{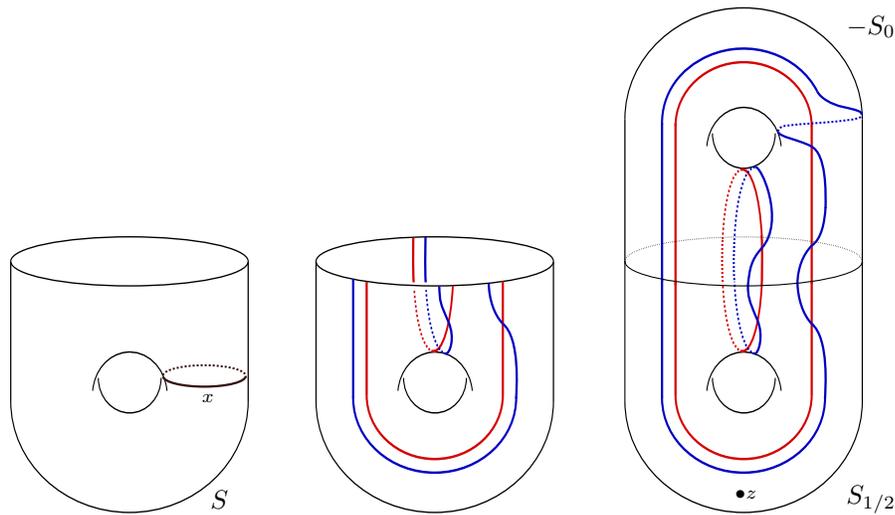}
\label{fig:arcs}
\caption{\quad On the left is the surface $S=S_{1,1}$. The figure in the middle shows the arcs $a_i$ (in red) and $b_i$ (in blue). On the right is a standard pointed Heegaard diagram for the open book $(S, D_x)$, where $D_x$ is a right-handed Dehn twist around the curve $x$.}
\end{center}
\end{figure}

For each $i=1, \dots, n$, let $y_i$ be the intersection point on $S_{1/2}$ between $\alpha_i$ and $\beta_i$. Then $\by = \{y_1,\dots,y_n\}$ represents an intersection point between $\Tb$ and $\Ta$ in $Sym^n(\Sigma)$, and we may think of $[\by,0]$ as an element of $\cfp(\Sigma,\beta,\alpha,z) = \cfp(-M_{S,\phi})$. 

\begin{theorem}[{\rm \cite[Theorem 3.1]{hkm3}}]
\label{thm:contactclass}
The image of $[\by,0]$ in $\hfp(-M_{S,\phi})$ is the Ozsv{\'a}th-Szab{\'o} contact class $c^+(S,\phi)$.
\end{theorem}

Now suppose that $S=S_{g,r}$, and let $B$ denote the boundary component of $S$ that we wish to cap off. Let $a_1,\dots, a_n$ (where $n=2g+r-1$) be pairwise disjoint, properly embedded arcs on $S$ so that $S \setminus \cup a_i$ is a disk and only $a_1$ intersects $B$. For each $i=1,\dots, n$, let $b_i$ be an arc obtained by changing $a_i$ via a small isotopy as described above. For each $i=2,\dots,n$, let $c_i$ be an arc obtained by changing $b_i$ via a similar isotopy (so that $c_i$ intersects each of $a_i$ and $b_i$ transversely in one point and with positive sign), and let $c_1$ be a curve on $S$ parallel to the boundary component $B$. See Figure \ref{fig:page} for an illustration of the curve $c_1$ and the arcs $a_i$.

For $i=1,\dots,n$, let $\alpha_i$ and $\beta_i$ be the curves on $\Sigma = S_{1/2}\cup -S_0$ defined by $$\alpha_i = a_i \times \{1/2\} \cup a_i \times \{0\},$$ $$\beta_i = b_i \times \{1/2\} \cup \phi(b_i) \times \{0\}$$ as above. In addition, define $$\gamma_1 = c_1 \times \{1/2\},$$ and let $$\gamma_i = c_i \times \{1/2\} \cup \phi(c_i) \times \{0\}$$ for $i =2,\dots,n$. Finally, place a basepoint $z$ in the ``big" region of $S_{1/2}\setminus\cup \alpha_i \setminus\cup \beta_i \setminus\cup \gamma_i$ (that is, neither in one of the thin strip regions nor in the region between $B$ and $\gamma_1$), and let $\alpha$, $\beta$ and $\gamma$ denote the sets of attaching curves $\{\alpha_1,\dots,\alpha_n\}$, $\{\beta_1,\dots,\beta_n\}$ and $\{\gamma_1,\dots,\gamma_n\}$. Then $(\Sigma, \alpha, \beta,z)$ is a standard pointed Heegaard diagram for $(S_{g,r},\phi)$. 

Let $K_B$ denote the binding component in $M_{S_{g,r},\phi}$ which corresponds to $B$. Observe that $\beta_1$ is a meridian of $K_B$, and that the Heegaard diagram $(\Sigma,\alpha, \beta \setminus\beta_1)$ specifies the knot complement $M_{S_{g,r},\phi} \setminus K_B$. Since $\gamma_1$ is a $0$-framed longitude of $K_B$ and $\gamma_i$ is isotopic to $\beta_i$ for $i\geq 2$, it follows that $(\Sigma,\alpha,\gamma)$ is a Heegaard diagram for the 3-manifold $M_{S_{g,r-1},\phi'}$ obtained by performing $0$-surgery on $K_B$. In fact, it is easy to see that $(\Sigma, \alpha, \gamma, z)$ is the stabilization of a standard pointed Heegaard diagram for the open book $(S_{g,r-1},\phi')$. 

For $i =1,\dots, n$, let $\theta_i$, $x_i$ and $y_i$ be the points in $\Sigma$ defined by $$\theta_i = \beta_i \cap \gamma_i \cap S_{1/2},$$ $$x_i = \gamma_i \cap \alpha_i \cap S_{1/2},$$ $$y_i = \beta_i \cap \alpha_i \cap S_{1/2},$$ and let $\bt$, $\mathbf{x}$ and $\mathbf{y}$ be the corresponding points in $Sym^n(\Sigma)$ defined by $$\mathbf{\Theta} = \{\theta_1,\dots,\theta_n\} \in \mathbb{T}_{\beta} \cap \mathbb{T}_{\gamma},$$ $$\mathbf{x} = \{x_1,\dots,x_n\} \in \mathbb{T}_{\gamma} \cap \mathbb{T}_{\alpha},$$ $$\mathbf{y} = \{y_1,\dots,y_n\} \in \mathbb{T}_{\beta} \cap \mathbb{T}_{\alpha}.$$ According to Theorem \ref{thm:contactclass}, the image of $[\mathbf{y},0]$ in $$\hfp(\Sigma, \beta, \alpha,z) = \hfp(-M_{S_{g,r},\phi})$$ is the contact class $c^+(S_{g,r},\phi)$; likewise, the image of $[\mathbf{x},0]$ in $$\hfp(\Sigma,\gamma, \alpha, z) = \hfp(-M_{S_{g,r-1},\phi'})$$ is $c^+(S_{g,r-1},\phi')$. Meanwhile, $[\bt,0]$ represents the top-dimensional generator of $$HF^{\leq0}(\Sigma,\beta,\gamma,z) = HF^{\leq 0}(\#^{n-1}(S^1 \times S^2)).$$

\subsection{The map induced by capping off}
\label{ssec:cap}
Suppose that $W$ is the cobordism from $M_{S_{g,r},\phi}$ to $M_{S_{g,r-1},\phi'}$ obtained by attaching a $0$-framed 2-handle to the knot $K_B$ in $M_{S_{g,r},\phi}$. As mentioned in the introduction, $W$ may be viewed as a cobordism from $-M_{S_{g,r-1},\phi'}$ to $-M_{S_{g,r},\phi}$ instead. If $\spc$ is a $\Sc$ structure on $W$, then the map $$F^+_{W,\spc}:\hfp(-M_{S_{g,r-1},\phi'})\rightarrow \hfp(-M_{S_{g,r},\phi})$$ is induced by the chain map, $$f^+_{W,\spc}:\cfp(-M_{S_{g,r-1},\phi'})\rightarrow \cfp(-M_{S_{g,r},\phi}),$$ which is defined using the pointed triple-diagram $(\Sigma, \beta, \gamma, \alpha, z)$. (Technically, this is a \emph{left-subordinate} triple diagram for the cobordism $W$, as opposed to the more often used notion of a \emph{right-subordinate} triple-diagram. Right- and left-subordinate diagrams induce the same maps on homology \cite[Lemma 5.2]{osz5}.) Recall that, for $\mathbf{v} \in \mathbb{T}_{\gamma} \cap \mathbb{T}_{\alpha}$, 
\begin{equation}
\label{eqn:map}f^+_{W,\spc} ([\mathbf{v},i]) = \sum_{\mathbf{w} \in  \mathbb{T}_{\beta} \cap \mathbb{T}_{\alpha}} \,\sum_{\{\psi \in \pi_2(\bt, \mathbf{v}, \mathbf{w})\,|\, \mu(\psi)=0, \,\spc_z(\psi) = \spc\}} (\#\mathcal{M}(\psi))\cdot[\mathbf{w},i-n_z(\psi)].
\end{equation} In this sum, $\pi_2(\bt, \mathbf{v}, \mathbf{w})$ is the set of homotopy classes of Whitney triangles connecting $\bt$, $\mathbf{v}$, and $\mathbf{w}$; $\mu(\psi)$ is the expected dimension of the moduli space, $\mathcal{M}(\psi)$, of holomorphic representatives of $\psi$; $\spc_z(\psi)$ is the $\Sc$ structure on $W$ corresponding to $\psi$; and $n_z(\psi)$ is the algebraic intersection number of $\psi$ with the subvariety $\{z\} \times Sym^{n-1}(\Sigma) \subset Sym^n(\Sigma)$. Below, we review some relevant definitions; for more details, see \cite{osz8}.

\begin{figure}[!htbp]
\labellist 
\hair 2pt 
\small\pinlabel $\mathbf{v}$ at 354 9
\pinlabel $\mathbf{w}$ at -16 9 
\pinlabel $\mathbf{r}$ at 175 330 
\pinlabel $\mathbb{T}_{\alpha}$ at 175 -6
\pinlabel $\mathbb{T}_{\beta}$ at 55 173 
\pinlabel $\mathbb{T}_{\gamma}$ at 282 173 
\pinlabel $u(\Delta)$ at 171 115
\endlabellist 
\begin{center}
\includegraphics[width=3.4cm]{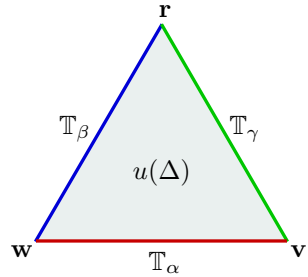}
\caption{\quad A Whitney triangle.}
\label{fig:whitney}
\end{center}
\end{figure}

Let $\Delta$ denote the 2-simplex with vertices $v_{\beta}$, $v_{\gamma}$ and $v_{\alpha}$ labeled clockwise, and let $e_{\beta}$, $e_{\gamma}$ and $e_{\alpha}$, respectively, denote the edges opposite these vertices. A \emph{Whitney triangle} connecting points $\mathbf{r},$ $ \mathbf{v}$ and $\mathbf{w}$ in $\Tb\cap \Tc$, $\Tc\cap\Ta$ and $\Tb \cap \Ta$ is a smooth map $$u:\Delta \rightarrow Sym^n(\Sigma)$$ with the boundary conditions that $u(v_{\alpha}) = \mathbf{r}$, $u(v_{\beta}) = \mathbf{v}$ and $u(v_{\gamma}) = \mathbf{w}$, and $u(e_{\beta}) \subset \Tb$, $u(e_{\gamma}) \subset \Tc$ and $u(e_{\alpha}) \subset \Ta$. See Figure \ref{fig:whitney} for a schematic depiction of this map.

Let $\mathcal{D}_1, \dots, \mathcal{D}_k$ denote the connected regions of $\Sigma \setminus \cup \alpha_i \setminus \cup \beta_i \setminus \cup \gamma_i$. A \emph{triply-periodic domain} for the pointed Heegaard diagram $(\Sigma, \beta, \gamma, \alpha, z)$ is a 2-chain $\mathcal{P} = \sum_i a_i \mathcal{D}_i$ in $C(\Sigma; \mathbb{Z})$ whose boundary is a sum of $\beta$, $\gamma$ and $\alpha$ curves, and whose multiplicity at the basepoint $z$ is 0 (the \emph{multiplicity} of a 2-chain at a point refers to the coefficient, in the 2-chain, of the region containing that point). The diagram $(\Sigma, \beta, \gamma, \alpha, z)$ is said to be \emph{weakly-admissible} if every non-trivial triply-periodic domain has both positive and negative multiplicities (this is slightly stronger than the definition of weak-admissibility given in \cite{osz8}). In general, the map $f^+_{W,\spc}$ is not well-defined unless the pointed triple-diagram $(\Sigma, \beta, \gamma, \alpha)$ is weakly-admissible. This is, therefore, our first consideration.

\begin{lemma}
\label{lem:weakadm}
The pointed triple-diagram $(\Sigma,\beta,\gamma,\alpha,z)$ constructed above is weakly-admissible.
\end{lemma}

\begin{proof}[Proof of Lemma \ref{lem:weakadm}]
Figure \ref{fig:page3} shows a local picture of $\Sigma$ near the intersection points $\theta_i$, $x_i$ and $y_i$ for $i \geq 2$. Let $\mathcal{P}$ be a triply-periodic domain whose multiplicities in the regions $A$, $B$, $C$, $D$, $E$ and $F$ are given by the integers $a$, $b$, $c$, $d$, $e$ and $f$, respectively. Note that $c = 0$ since the region $C$ contains the basepoint. Since $\partial \mathcal{P}$ consists of complete $\beta$, $\gamma$ and $\alpha$ curves, it must be that $$b= d-e = -f,$$ $$a = b-d = -e.$$ Therefore, $\mathcal{P}$ has both positive and negative multiplicities unless $$a = b=c=d=e=f=0.$$ We perform this local analysis for each $i=2,\dots,n$ and conclude that either $\mathcal{P}$ has both positive and negative multiplicities or $\partial \mathcal{P}$ is a linear combination of the curves $\beta_1$, $\gamma_1$ and $\alpha_1$. Let us assume the latter.

\begin{figure}[!htbp]
 
\labellist 
\hair 2pt 
\pinlabel $\alpha_i$ at 102 187 
\pinlabel $\beta_i$ at 66 187 
\pinlabel $\gamma_i$ at 24 187 

\pinlabel $S_{1/2}$ at 167 145
\pinlabel $-S_{0}$ at 164 187

\small\pinlabel $\bullet$ at 63 129
\small\pinlabel $\theta_i$ at 56 123
\pinlabel $\bullet z$ at 45 25 
\pinlabel $\bullet$ at 101 113
\pinlabel $x_i$ at 110 119
 
\pinlabel $\bullet$ at 100 67
\pinlabel $y_i$ at 92 60

\tiny
\pinlabel $A$ at 47 154
\pinlabel $B$ at 82 144
\pinlabel $C$ at 137 119
\pinlabel $C$ at 49 79
\pinlabel $D$ at 86 100
\pinlabel $E$ at 123 78
\pinlabel $F$ at 119 25


\endlabellist 

\begin{center}
\includegraphics[width=5.2cm]{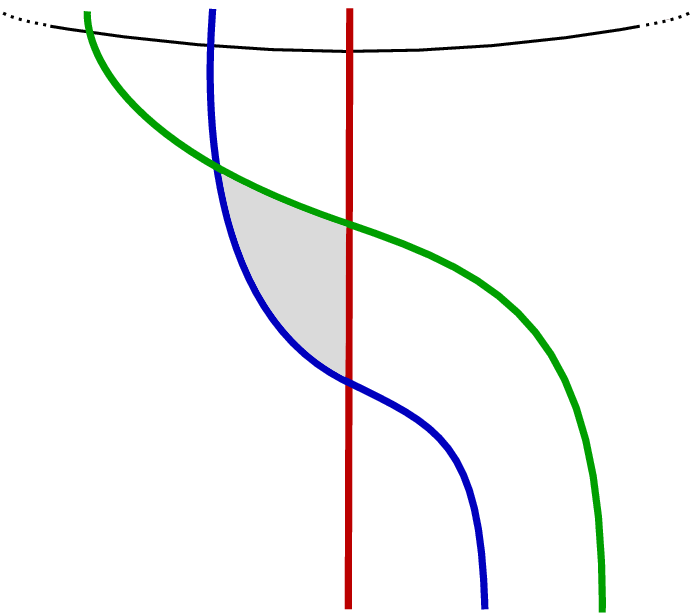}
\caption{\quad The local picture near the intersection points $\theta_i$, $x_i$ and $y_i$. The shaded region is $\Delta_i$.}
\label{fig:page3}
\end{center}
\end{figure}

Now, consider the regions labeled $A$, $B$, and $C$ in Figure \ref{fig:page9}, and suppose that $\mathcal{P}$ has multiplicities $a$, $b$ and $c$ in these regions. Again, $c=0$ since region $C$ contains the basepoint; and $a=-b$. Therefore, $\mathcal{P}$ has both positive and negative multiplicities unless $a=b=0$, in which case $\partial \mathcal{P}$ is some multiple of the curve $\gamma_1$. Since $\gamma_1$ is not null-homologous in $\Sigma$, this multiple must be zero, which implies that $\mathcal{P}$ is the trivial domain. To summarize, we have shown that $\mathcal{P}$ has both positive and negative multiplicities unless $\mathcal{P}$ is trivial. Hence, the diagram $(\Sigma, \beta,\gamma,\alpha,z)$ is weakly-admissible.

\begin{figure}[!htbp]

\labellist 
\hair 2pt 
\pinlabel $\alpha_1$ at 185 217
\pinlabel $\beta_1$ at 68 180 
\pinlabel $\gamma_1$ at -10 117 

\pinlabel $S_{1/2}$ at -6 40 
\pinlabel $-S_0$ at -9 194

\small\pinlabel $\bullet$ at 81 92
\pinlabel $\theta_1$ at 77 84
\pinlabel $\bullet z$ at 60 25 

\pinlabel $\bullet$ at 107 96
\pinlabel $x_1$ at 118  91
 
\pinlabel $\bullet$ at 253 96
\pinlabel $y_1$ at 263 98
\tiny
\pinlabel $A$ at 185 29
\pinlabel $B$ at 250.5 117
\pinlabel $C$ at 270 55
\pinlabel $C$ at 224 76
\pinlabel $D$ at 92 106
\pinlabel $E$ at 67 102
\pinlabel $E$ at 115 109

\endlabellist 

\begin{center}
\includegraphics[width=7.5cm]{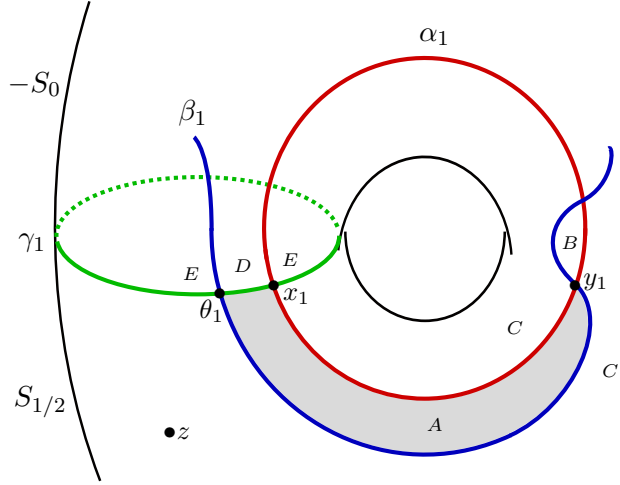}
\caption{\quad The local picture near the intersection points $\theta_1$, $x_1$ and $y_1$. The shaded region is $\Delta_1$.}
\label{fig:page9}
\end{center}
\end{figure}
\end{proof}

Recall that a homotopy class $\psi$ of Whitney triangles has an associated domain $\mathcal{D}(\psi) = \sum_i n_{p_i}(\psi) \mathcal{D}_i$, where $p_i$ is a point in $\mathcal{D}_i$. For each $i=1,\dots,n$, let $\Delta_i \subset S_{1/2}$ be the shaded triangular region with vertices at $\theta_i$, $x_i$ and $y_i$ shown in Figures \ref{fig:page3} and \ref{fig:page9}. Then the homotopy class $\psi_0 \in \pi_2(\bt,\bx,\by)$ with domain $\mathcal{D}(\psi_0) = \Delta_1 + \dots + \Delta_n$ has a unique holomorphic representative, by the Riemann Mapping Theorem (in particular, $\mu(\psi_0) = 0$). Let $\spc_0$ denote the $\Sc$ structure $\spc_z(\psi_0)$, and observe that $n_z(\psi_0)=0$.

\begin{proposition}
Suppose that $\psi$ is a homotopy class of Whitney triangles connecting $\mathbf{\Theta},$ $\bx$ and any other point $\bw \in \Tb\cap\Ta$. Let $w_i$ denote the component of $\bw$ on $\gamma_i$. If $\psi$ has a holomorphic representative and satisfies $n_z(\psi)=0$, then $w_i = y_i$ for $i = 2,\dots, n$, and $\mathcal{D}(\psi) = \Delta_1'+\Delta_2+\dots+\Delta_n$, where $\Delta_1'$ is a (possibly non-embedded) triangle in $\Sigma \setminus \{z\}$ with vertices at $\theta_1$, $x_1$, and $w_1$. If, in addition, $\spc_z(\psi) = \spc_0$, then $\psi = \psi_0$ and $\bw = \by$.
\label{prop:nat}

\end{proposition}

This proposition implies that the map $f^+_{W,\spc_0}$ sends $[\bx,0]$ to $[\by,0]$, proving Theorem \ref{thm:nat}.

\begin{proof}[Proof of Proposition \ref{prop:nat}]
Suppose $\psi$ has a holomorphic representative and satisfies $n_z(\psi)=0$. Then every coefficient in the domain $\mathcal{D}(\psi)$ is non-negative, and $\mathcal{D}(\psi)$ must have multiplicity 0 in the region containing the basepoint $z$. Moreover, the oriented boundary of $\mathcal{D}(\psi)$ consists of arcs along the $\beta$ curves from the points $w_1,\dots,w_n$ to the points $\theta_1,\dots,\theta_n$; arcs along the $\gamma$ curves from the points $\theta_1,\dots,\theta_n$ to the points $x_1,\dots,x_n$; and arcs along the $\alpha$ curves from the points $x_1,\dots,x_n$ to the points $w_1,\dots,w_n$.

Let $a$, $b$, $c$, $d$, $e$ and $f$ be the multiplicities of $\mathcal{D}(\psi)$ in the regions $A$, $B$, $C$, $D$, $E$ and $F$ shown in Figure \ref{fig:page3}. We have already established that $c=0$. The boundary constraints on $\mathcal{D}$ then imply that \begin{equation}\label{eqn:one} a+d = b+1,\end{equation}  $$d = b+e+1.$$ Subtracting one equation from the other, we find that $a = -e$. Since all coefficients of $\mathcal{D}(\psi)$ are non-negative, $a=e=0$. If $w_i \neq y_i,$ then the constraints on $\partial \mathcal{D}(\psi)$ force $f+d =0$, which implies that $f=d=0$. However, plugging this back into Equation \ref{eqn:one}, together with $a=0,$ implies that $0=b+1$, which contradicts the fact that $b$ is non-negative. As a result, it must be the case that $w_i=y_i$. Then the constraints on $\partial \mathcal{D}(\psi)$ (together with the fact that $e=c=0$) require that $d+f=1.$ Combined with Equation \ref{eqn:one}, this implies that $d=1$ and $f=b=0$. So, we have found that $d=1$ and $a=b=c=e=f=0$; that is, the domain $\mathcal{D}(\psi)$ is locally just $\Delta_i$.

We perform this local analysis for each $i = 2,\dots,n$ and conclude that $w_i = y_i$ for $i=2,\dots,n$ and that $\mathcal{D}(\psi) = \Delta_1' + \Delta_2 + \dots + \Delta_n$, where $\Delta_1'$ is a region whose oriented boundary consists of arcs along $\beta_1$ from $w_1$ to $\theta_1$; along $\gamma_1$ from $\theta_1$ to $x_1$; and along $\alpha_1$ from $x_1$ to $w_1$. In fact, since $\Delta_2, \dots, \Delta_n$ are triangles in $\Sigma$ and $\mathcal{D}(\psi)$ is the image of a map from the $n$-fold branched cover of a triangle into $\Sigma$ (see \cite{osz8}), $\Delta_1'$ must be a (possibly non-embedded) triangle in $\Sigma$ as well which avoids the basepoint $z$.

Now, suppose that $\psi$ also satisfies $\spc_z(\psi) = \spc_0$. The only thing left to prove is that $w_1 = y_1$ and $\Delta_1' = \Delta_1$. Since $\spc_z(\psi) = \spc_0=\spc_z(\psi_0)$, it follows from \cite[Proposition 8.5]{osz8} that $$\mathcal{D}(\psi)-\mathcal{D}(\psi_0) = \mathcal{D}(\phi_1) + \mathcal{D}(\phi_2) + \mathcal{D}(\phi_3),$$ where $\phi_1$, $\phi_2$ and $\phi_3$ are homotopy classes of Whitney disks in $\pi_2(\by,\bw)$, $\pi_2(\bt,\bt)$ and $\pi_2(\bx,\bx)$, respectively. Since $\mathcal{D}(\psi)-\mathcal{D}(\psi_0) = \Delta_1'-\Delta_1$, and $\gamma_1$ is homologically independent of both $\alpha_1$ and $\beta_1$, $\mathcal{D}(\phi_2) = \mathcal{D}(\phi_3)=0$, and $\mathcal{D}(\phi_1)$ is a disk in $\Sigma\setminus \{z\}$ whose oriented boundary consists of arcs along $\alpha_1$ from $y_1$ to $w_1$, and arcs along $\beta_1$ from $w_1$ to $y_1$. 

Let $a$, $b$, $d$ and $e$ be the multiplicities of $\Delta_1'$ in the regions labeled $A$, $B$, $D$ and $E$ in Figure \ref{fig:page9} (the multiplicity of $\Delta_1'$ in region $C$ is 0). Since $\Delta_1'-\Delta$ is the disk $\mathcal{D}(\phi_1)$, the multiplicity of $\Delta_1'-\Delta$ in the region labeled $D$ must be the same as its multiplicity in the region labeled $A$; that is, \begin{equation}\label{eqn:two}d=a-1\end{equation} (the multiplicities of $\Delta$ in these regions are 0 and 1, respectively). But the boundary constraints on $\mathcal{D}(\psi)$ imply that $$a+e = d+1.$$ Combining this equation with the former, we find that $e=0$. If $w_1 \neq y_1$, then the same boundary constraints require that $a+b=0$. Yet, combined with Equation \ref{eqn:two}, this implies that either $a$ or $b$ is negative, which contradicts our assumption that $\psi$ has a holomorphic representative. Therefore, $w_1=y_1$, and the boundary constraints on $\mathcal{D}(\psi)$ imply that $$a+b=1.$$ It follows that $a=1$ and $b=d=0$, and, hence, that $\Delta_1'=\Delta$. Thus, $\psi = \psi_0$, completing the proof of Proposition \ref{prop:nat}.
 \end{proof}

\section{Contact surgery and Legendrian stabilization}
\label{sec:legstab}
\begin{figure}[!htbp]

\labellist 
\hair 2pt 
\small\pinlabel $K$ at 60 2
\pinlabel $S_+(K)$ at 240 2
\pinlabel $S_-(K)$ at 420 2

\endlabellist 

\begin{center}
\includegraphics[width=8cm]{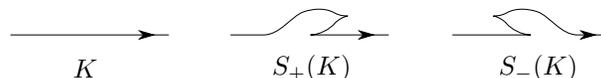}
\caption{\quad The Legendrian stabilizations $S_+(K)$ and $S_-(K),$ as seen via their front projections in $N$.}
\label{fig:stab2}
\end{center}
\end{figure}

In this section, we describe how contact $\pm 1$-surgery on a stabilized Legendrian knot fits into the framework of capping off. Suppose that $K$ is an oriented Legendrian knot in a contact 3-manifold $(M,\xi)$, and let $\mathbb{R}^3_y$ denote the quotient $\mathbb{R}^3/(y \sim y+1).$ There is a contactomorphism from a neighborhood of $K$ to $(N,\xi')$, where $$N = \{(x,y,z)\in \mathbb{R}^3_y\, | \, x^2+z^2<\epsilon\},$$ $\xi' = \text{ker} (dz+xdy),$ and $K$ is sent to the image of the $y$-axis in $N$. There is a natural ``front" projection in $N$ defined by the map which sends $(x,y,z)$ to $(y,z)$. In \cite{EH2}, Etnyre and Honda define the positive and negative Legendrian stabilizations of $K$, $S_+(K)$ and $S_-(K)$, to be the Legendrian knots in $M$ corresponding to the curves in $N$ shown in Figure \ref{fig:stab2}. Note that this definition agrees with the usual definition of stabilization for Legendrian knots in the standard tight contact structure on $S^3$.

By incorporating $K$ into the 1-skeleton of a contact cell decomposition for $(M,\xi)$, we can find an open book $(S_{g,r},\phi)$ compatible with $\xi$ so that $K$ lies on a page of this open book and the contact framing of $K$ agrees with the framing induced by this page. The lemma below is based upon this idea as well.

\begin{figure}[!htbp]

\labellist 
\hair 2pt 
\tiny\pinlabel $K$ at -12 92
\pinlabel $S_+(K)$ at 272 102
\pinlabel $S_-(K)$ at 272 79
\pinlabel $B_+$ at 410 128
\pinlabel $B_-$ at 410 54

\endlabellist 

\begin{center}
\includegraphics[width=10.5cm]{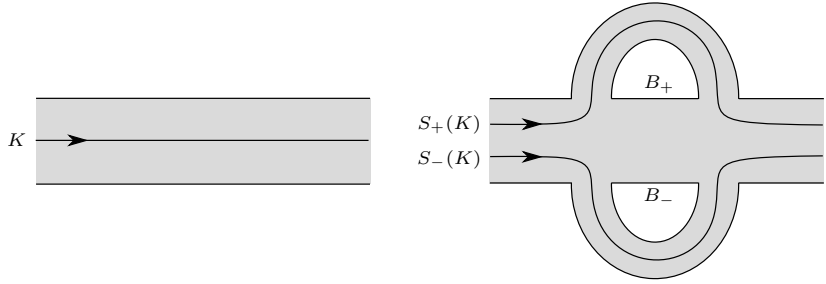}
\caption{\quad On the left is a neighborhood of a piece of $K$ in $S_{g,r}$. On the right is a portion of the twice stabilized open book with the curves $S_+(K)$ and $S_-(K)$. We have labeled the two new boundary components $B_+$ and $B_-$. }
\label{fig:stabpage}
\end{center}
\end{figure}

\begin{lemma}[{\rm \cite[Lemma 3.3]{et}}]
\label{lem:stab}
Suppose the oriented Legendrian knot $K$ lies on a page of the open book $(S_{g,r},\phi)$. If we positively stabilize $(S_{g,r},\phi)$ twice as shown in Figure \ref{fig:stabpage}, then we may isotop the page of the stabilized open book so that both $S_+(K)$ and $S_-(K)$ appear on the page as in Figure \ref{fig:stabpage}. The contact framings of these stabilized Legendrian knots agree with their page framings.
\end{lemma}

\begin{proof}[Proof of Theorem \ref{thm:stab}]
Let $(S_{g,r+2}, \phi')$ refer to the twice stabilized open book in Lemma \ref{lem:stab}, and let $K'$ be either $S_+(K)$ or $S_-(K)$. We think of $K'$ as lying in a page of this open book, per Lemma \ref{lem:stab}. Since the contact framings of $K$ and $K'$ agree with their page framings, the contact 3-manifolds $(M_{\pm1}(K), \xi_{\pm 1}(K))$ and $(M_{\pm1}(K'), \xi_{\pm 1}(K'))$ are supported by the open books $(S_{g,r},\phi \cdot t_{K}^{\mp 1})$ and $(S_{g,r+2},\phi' \cdot t_{K'}^{\mp 1})$, respectively. Note that $(S_{g,r},\phi \cdot t_{K}^{\mp 1})$ is obtained from $(S_{g,r+2},\phi' \cdot t_{K'}^{\mp 1})$ by capping off the boundary components $B_+$ and $B_-$. Therefore, by two applications of Theorem \ref{thm:nat}, there is a $U$-equivariant map $$F^+:\hfp(-M_{\pm1}(K))\rightarrow \hfp(-M_{\pm 1}(K'))$$ which sends $c^+(\xi_{\pm 1}(K))$ to $c^+(\xi_{\pm 1}(K'))$.
\end{proof}

Theorem \ref{thm:stab} has a nice interpretation in terms of contact surgery. Recall that, for $n\in \mathbb{Z}^{<0}$, contact $n$-surgery on a Legendrian knot $K\subset M$ may be performed by stabilizing $K$ a total of $-n-1$ times to obtain $K'$, and then performing contact $-1$-surgery on $K'$ \cite{DG2,DGS}. In particular, such contact surgery is not unique unless $n=-1$; the ambiguity corresponds to the various ways of stabilizing $K$. By applying Theorem \ref{thm:stab} to $K'$, we obtain, under the appropriate interpretations of the contact manifolds $(M_{n}(K),\xi_n(K))$ and $(M_{n-1}(K),\xi_{n-1}(K))$, a map $$F^+:\hfp(-M_{n}(K))\rightarrow \hfp(-M_{n-1}(K))$$ which sends $c^+(\xi_{n}(K))$ to $c^+(\xi_{n-1}(K))$.


\section{Gluing open books}
\label{sec:glue}

\begin{figure}[!htbp]

\labellist 
\hair 2pt 
\small
\pinlabel $B$ at 205 240
\pinlabel $B'$ at 207 147
\pinlabel $B\#B''$ at 515 213

\endlabellist 

\begin{center}
\includegraphics[width=8.5cm]{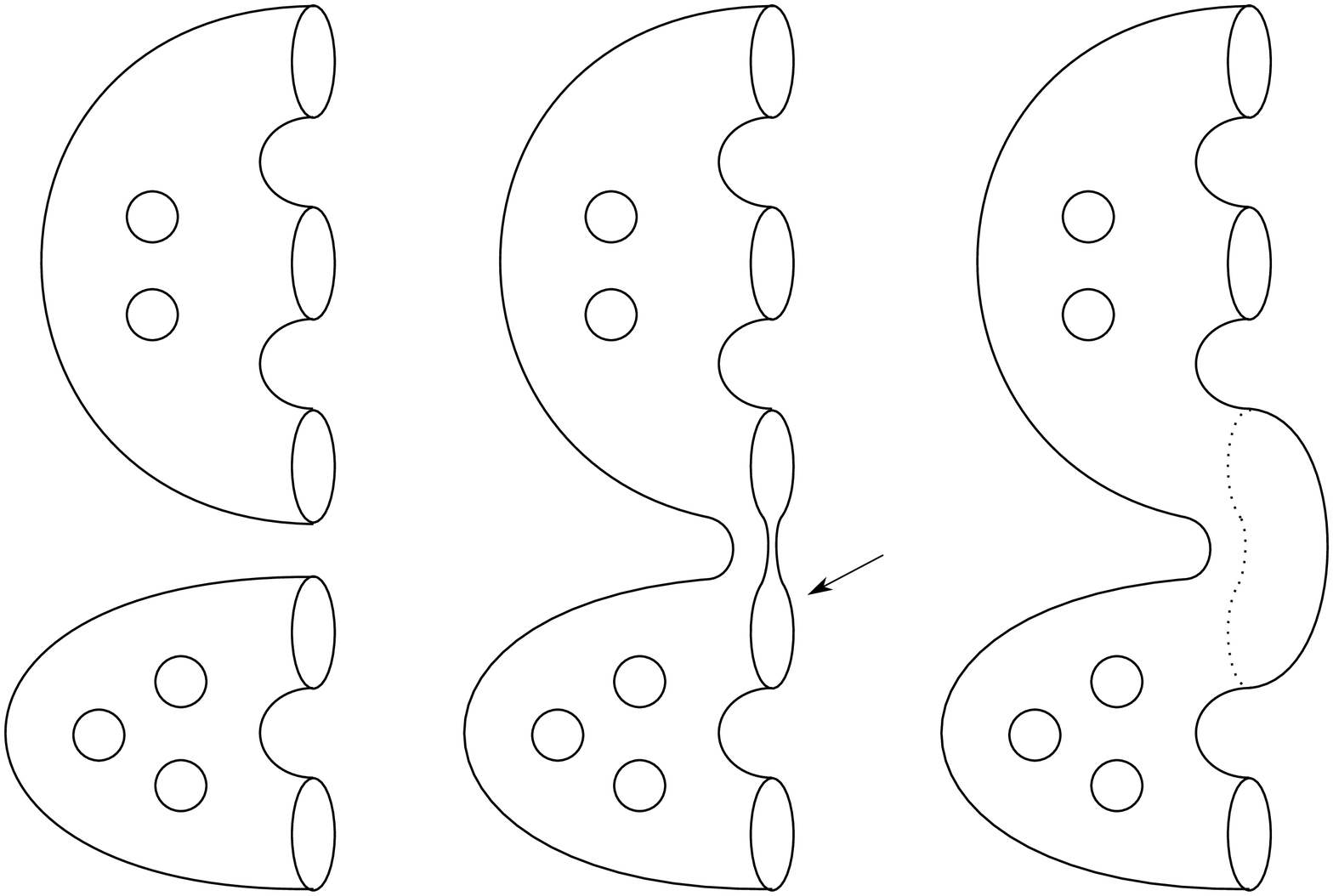}
\caption{\quad An example of gluing via boundary connected sum and capping off.}
\label{fig:glue}
\end{center}
\end{figure}

Let $(S'',\phi'')$ denote the result of gluing $(S,\phi)$ to $(S',\phi')$ along boundary components $B$ and $B'$ of $S$ and $S'$, respectively. The open book $(S'',\phi'')$ may also be obtained by taking the boundary connected sum of the open books $(S,\phi)$ and $(S',\phi')$ along $B$ and $B'$, and then capping off the boundary component $B\,\#\,B'$ of the resulting surface, as illustrated in Figure \ref{fig:glue}. Since the boundary connected sum of these open books supports the contact connected sum $\xi_{S,\phi}\,\#\,\xi_{S',\phi'}$, Theorem \ref{thm:nat} and \cite[Proposition 2.1]{osz1} imply the following.

\begin{lemma}
\label{lem:glue}
Suppose that $(S'',\phi'')$ is the open book obtained by gluing $(S,\phi)$ to $(S',\phi')$ along boundary components $B, B'$. If $c(S,\phi)$ and $c(S',\phi')$ are both non-zero, then so is $c(S'',\phi'')$.
\end{lemma}

\begin{figure}[!htbp]

\labellist 
\hair 2pt 
\small

\pinlabel $B$ at 123 170
\pinlabel $B'$ at 122 60
\pinlabel $B\#B''$ at 445 23
\endlabellist 

\begin{center}
\includegraphics[width=8.5cm]{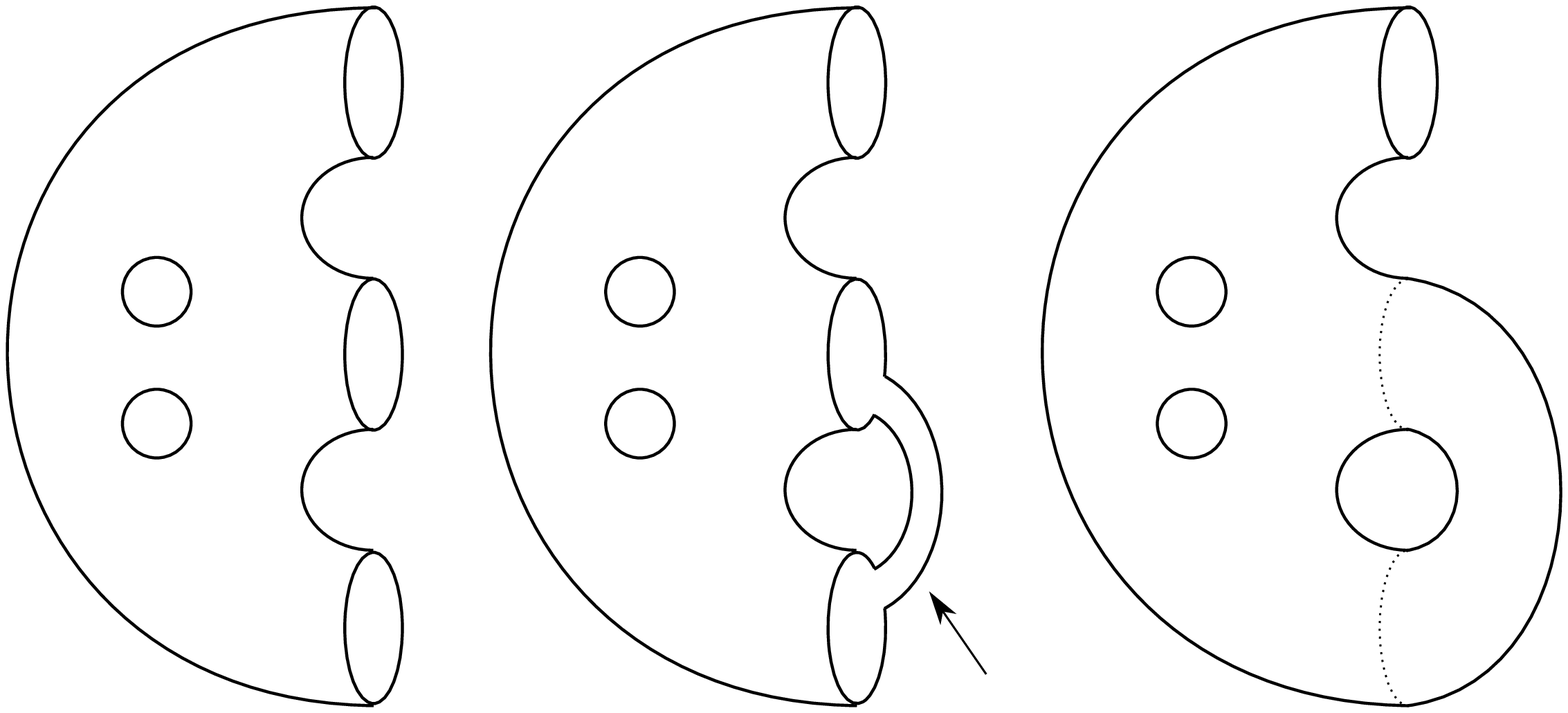}
\caption{\quad An example of self-gluing via 1-handle attachment and capping off.}
\label{fig:1hsum}
\end{center}
\end{figure}

One further observation is needed to complete the proof of Theorem \ref{thm:glue}. Namely, suppose that $(S,\phi)$ is an open book with at least three boundary components. Let $B$ and $B'$ denote two of them, and consider the open book $(S',\phi')$ obtained from $(S,\phi)$ by gluing $B$ to $B'$ (we shall refer to this operation as \emph{self-gluing}). $(S',\phi')$ may alternatively be obtained by attaching a 1-handle to $(S,\phi)$ with feet on $B$ and $B'$, and then capping off the boundary component $B\,\#\,B'$ of the resulting surface; see Figure \ref{fig:1hsum}. Since this 1-handle attachment corresponds to taking a contact connected sum with the Stein fillable contact structure on $S^1\times S^2$, Theorem \ref{thm:nat} and \cite[Proposition 2.1]{osz1} combine to give lemma below.

\begin{lemma}
\label{lem:sglue}
Suppose that $(S',\phi')$ is the open book obtained from $(S,\phi)$ by self-gluing along $B,B'$. If $c(S,\phi)$ and $c(S',\phi')$ are both non-zero, then so is $c(S'',\phi'')$.
\end{lemma}

Now, suppose that the open book $(S'',\phi'')$ is obtained by gluing $(S,\phi)$ to $(S',\phi')$ along boundary components $B_1,\dots,B_n$ of $S$ and $B'_1,\dots,B'_n$ of $S'$, as in the introduction. Note that $(S'',\phi'')$ is result of gluing $(S,\phi)$ to $(S',\phi')$ along $B_1,B_1'$, followed by $n-1$ self-gluings along the other $B_i,B'_i$. Theorem \ref{thm:glue} therefore follows from Lemmas \ref{lem:glue} and \ref{lem:sglue}.

\begin{remark} Gluing has an inverse operation called \emph{splitting}. More precisely, suppose that $\phi$ is a reducible diffeomorphism of $S$ which fixes disjoint simple closed curves $C_1,\dots,C_n$ pointwise. Splitting $S$ along the $C_i$, one obtains open books $(S^{(1)},\phi^{(1)}),\dots, (S^{(m)},\phi^{(m)})$; conversely, we can recover $(S,\phi)$ from the $(S^{(j)},\phi^{(j)})$ via a combination of gluings and self-gluings. Lemmas \ref{lem:glue} and \ref{lem:sglue} then tell us that $c(S,\phi)$ is non-zero as long as all of the $c(S^{(j)},\phi^{(j)})$ are.
\end{remark}

\section{The support genera of contact structures compatible with $(S_{1,1},\phi)$.}
\label{sec:sg}
The mapping class group of $S_{1,1}$ is generated by Dehn twists around the curves $a$ and $b$ shown on the right in Figure \ref{fig:example}. It is well-known that this group is isomorphic to the braid group $B_3$ by an isomorphism $\Phi: \text{MCG}(S_{1,1},\partial S_{1,1})\rightarrow B_3$ which sends the Dehn twists $t_a$ and $t_b$ to the standard generators $\sigma_1$ and $\sigma_2$ of $B_3$. So, by a theorem of Murasugi on 3-braids \cite{mur}, we have the following.

\begin{lemma}
\label{lem:pamaps}
Let $h = (t_at_b)^3$. Any diffeomorphism of $S_{1,1}$ which fixes the boundary pointwise and is freely isotopic to a pseudo-Anosov map is, up to conjugation, isotopic (rel. $\partial$) to a diffeomorphism $$\phi_{\bn,d} = h^d\cdot t_b^{}t_a^{-n_1}\cdots t_b^{}t_a^{-n_k}$$ for some $k$-tuple of non-negative integers $\bn=(n_1,\dots,n_k)$ for which some $n_i\neq 0$, and some $d\in\mathbb{Z}$.
\end{lemma}
 
The diffeomorphism $h$ represents a ``half-twist" around a curve $\delta$ parallel to the boundary of $S_{1,1}$; that is, $h^2 = t_{\delta}$. Let $\xi_{\bn,d}$ denote the contact structure compatible with the open book $(S_{1,1},\phi_{\bn,d}).$ In this short section, we prove the following.

\begin{proposition}
\label{prop:sg}
The support genus of $\xi_{\bn,d}$ is zero for $d\leq 0$, and one for $d>1$.
\end{proposition}

Note that this proposition is inconclusive for $d=1$.

\begin{proof}[Proof of Proposition \ref{prop:sg}]

In \cite{bald1, hkm3}, it is shown that $\xi_{\bn,d}$ is tight if and only if $d>0$. Recall from the introduction that $sg(\xi) = 0$ if $\xi$ is overtwisted \cite{et}. It therefore follows immediately that $sg(\xi_{\bn,d}) = 0$ for $d\leq 0$.

To simplify notation in this proof, we let $M_{\bn,d}$ denote the 3-manifold with open book decomposition $(S_{1,1},\phi_{\bn,d}).$ Observe that $M_{\bn,d}$ is the double cover of $S^3$ branched along the closed braid $B_{\bn,d}$ corresponding to the braid word $\Phi(\phi_{\bn,d}) \in B_3$ (see \cite[Section 2]{bald4}, for example). Note that $B_{\bn,d}$ is obtained from the alternating braid $B_{\bn,0}$ by adding $d$ full positive twists. It is clear that, as manifolds, links and contact structures, $M_{\bn,d}$, $B_{\bn,d}$ and $\xi_{\bn,d}$ are invariant under the action of cyclic permutation on the tuple $\bn$. 

The contact structure $\xi'$ in Example \ref{ex:sg} associated to the open book $(S_{1,1}, (t_at_b)^5 t_b^2)$ is simply $\xi_{(1),2}$ in our notation. In that example, we showed that there is some $d\in \mathbb{N}$ for which $c^+(\xi_{(1),2}) \notin U^d\cdot\hfp(-M_{(1),2})$; let us call this condition on $c^+(\xi_{(1),2})$ \emph{Condition} $\bR$. Recall that if $c^+(\xi)$ satisfies Condition $\bR$, then $sg(\xi)>0$, by Proposition \ref{prop:nontor}. 

For a $k$-tuple $\bn$ as above, let $\bn^-$ denote the $k$-tuple obtained from $\bn$ by adding 1 to its kth entry, and let $\bn^0$ denote the $(k+1)$-tuple obtained by concatenating $\bn$ with a 0. Starting from the $1$-tuple $\bn=(1)$, we can obtain any $k$-tuple of the form described in Lemma \ref{lem:pamaps} by repeated applications of the operations $\bn \mapsto \bn^-$, $\bn \mapsto \bn^0$, together with cyclic permutation. Moreover, for $d\geq2$, the monodromy $\phi_{\bn,d}$ is obtained from $ \phi_{\bn,2}$ by composition with $6d-12$ right handed Dehn twists around the curves $a$ and $b$. So, by the naturality of the contact invariant under maps induced by Stein cobordisms \cite{osz1,hkm3}, $c^+(\xi_{\bn,d})$ satisfies Condition $\bR$ as long as $c^+(\xi_{\bn,2})$ does.

Thus, in order to prove Proposition \ref{prop:sg}, it suffices to show that if $\xi_{\bn,2}$ satisfies Condition $\bR$, then so do $\xi_{\bn^-,2}$ and $\xi_{\bn^0,2}.$ For the latter, observe that $\phi_{\bn^0,2} = \phi_{\bn,2}\cdot t_b$. So, the naturality of the contact invariant under maps induced by Stein cobordisms implies that $c^+(\xi_{\bn^0,2})$ satisfies Condition $\bR$ as long as $c^+(\xi_{\bn,2})$ does. Proving the other implication takes slightly more work.

Observe that $\phi_{\bn,2} = \phi_{\bn^-,2}\cdot t_a$, and consider the map $\widehat{G}:\hf(-M_{\bn,2})\rightarrow \hf(-M_{\bn^-,2})$ induced by the corresponding Stein 2-handle cobordism (or, equivalently, by $-1$-surgery on a copy of the curve $a$ in the open book for $M_{\bn^-,2}$). To understand $\widehat{G}$, it helps to think of $M_{\bn,d}$ and $M_{\bn^-,d}$ as the branched double covers $\Sigma(B_{\bn,d})$ and $\Sigma(B_{\bn^-,d})$. Note that $B_{\bn,d}$ is obtained from $B_{\bn^-,d}$ by taking the \emph{oriented} resolution of $B_{\bn^-,d}$ at a negative crossing. Let us denote the \emph{unoriented} resolution at this crossing by $B^u_{\bn^- ,d}$. It is not hard to see that $B^u_{\bn^-,d}$ is an alternating link and does not depend on $d$. Moreover, the double covers of $S^3$ branched along these braids fit into the surgery exact triangle below (for $d=2$) \cite{osz12}.
$$\xymatrix @!0@C=5pc@R=5pc{
  \hf(-\Sigma(B_{\bn,2})) \ar[rr]^{\widehat{G}}&  &\hf(-\Sigma(B_{\bn^-,2}) )\ar[dl] \\
 &\hf(-\Sigma(B^u_{\bn^-, 2})). \ar[ul] }$$

Let $\bF$ denote the field with two elements, and let $\mathcal{T}^+$ denote the $\bF[U]$-module $\bF[U, U^{-1}]/\bF[U]$. From the grading calculations in \cite[Section 6]{bald4}, it follows that \begin{equation}\label{eqn:hfp}\hfp(-\Sigma(B_{\bn,2})) \cong (\mathcal{T}^+_0)^{\text{det}(B_{\bn,0})} \oplus \bF_1.\end{equation} Here, the subscripts denote absolute $\zzt$ gradings. The long exact sequence relating $\hf$ and $\hfp$ therefore implies that $$\hf(-\Sigma(B_{\bn,2})) \cong (\bF_0)^{\text{det}(B_{\bn,0})+1} \oplus \bF_1. $$ In particular, $$\text{rk}(\hf(-\Sigma(B_{\bn,2}))) = 2+\text{det}(B_{\bn,0}).$$ Of course, the analogous formula holds for $\text{rk}(\hf(-\Sigma(B_{\bn^-,2}) ))$. Moreover, since $B^u_{\bn^-,2} = B^u_{\bn^-,0}$ is alternating, we know from \cite{osz12} that $$\text{rk}(\hf(-\Sigma(B^u_{\bn^-,2})))= \text{det}(B^u_{\bn^-,0}).$$ And, since $B_{\bn^-,0}$ is an alternating link, its determinant satisfies \begin{equation}\label{eqn:detsum}\text{det}(B_{\bn^-,0}) = \text{det}(B_{\bn,0})+\text{det}(B^u_{\bn^-,0}).\end{equation} Combined with the rank formulae above, Equation \ref{eqn:detsum} implies that $$\text{rk}(\hf(-\Sigma(B_{\bn^-,2}) )) = \text{rk}(\hf(-\Sigma(B_{\bn,2}) ))+ \text{rk}(\hf(-\Sigma(B^u_{\bn^-,2}) )).$$ Therefore, our surgery exact triangle splits into a short exact sequence; in particular, the map $\widehat{G}$ is injective. 

Now, let us assume that $c^+(\xi_{\bn,2})$ satisfies Condition $\bR$. Then, according to Equation \ref{eqn:hfp}, $c^+(\xi_{\bn,2})$ must have absolute $\zzt$ grading 1 (and, hence, so does $c(\xi_{\bn,2})$). Since $\widehat{G}$ sends $c(\xi_{\bn,2})$ to $c(\xi_{\bn^-,2})$, and maps induced by cobordisms preserve relative $\zzt$ gradings (and are homogeneous with respect to these gradings) \cite{osz14}, the injectivity of $\widehat{G}$ forces $c(\xi_{\bn^-,2})$ (and, hence, $c^+(\xi_{\bn^-,2})$) to have absolute $\zzt$ grading 1 as well. But since $c^+(\xi_{\bn^-,2})$ has absolute $\zzt$ grading 1, it must satisfy Condition $\bR$, by the analogue of Equation \ref{eqn:hfp} for $\hfp(-\Sigma(B_{\bn^-,2}))$.
\end{proof}

\section{Capping off and periodic open books}
\label{sec:capper}

In this section, we study the 3-dimensional invariants associated to contact structures supported by genus one open books with periodic monodromy. 

\subsection{Periodic diffeomorphisms and the first Chern class} Recall that $\hfp(M,\spt)$ comes equipped with a $\mathbb{Q}$-grading whenever $c_1(\spt)$ is a torsion class. We denote the grading of a homogeneous element $x \in \hfp(M,\spt)$ by $gr(x)$. The proposition below appears in a slightly different form in \cite{osz1}.

\begin{proposition} [{\rm \cite[Proposition 4.6]{osz1}}]
\label{prop:grading} If $(M,\xi)$ is a contact 3-manifold for which $c_1(\spt_{\xi})$ is torsion, then $d_3(\xi) = -gr(c^+(\xi))-1/2.$
\end{proposition}

Suppose that $\phi$ is a diffeomorphism of $S_{g,r}$ such that $\phi^m$ is freely isotopic to the identity. Let $B_1,\dots, B_r$ denote the boundary components of $S_{g,r}$, and let $c_i$ be a curve on $S_{g,r}$ parallel to $B_i$ for each $i=1,\dots,r$. Since $\phi$ is periodic, $\phi^m$ is freely isotopic to the identity for some $m \in \mathbb{N}$. It follows that $\phi^m$ is isotopic to a product of Dehn twists of the form $t_{c_1}^{k_1}\cdots t_{c_r}^{k_r}$. For each $i=1,\dots,r$, we define the \emph{fractional Dehn twist coefficient} of $\phi$ around $B_i$ to be $k_i/m$ (see \cite{ch}). If $\phi$ is periodic, Colin and Honda show that the contact structure compatible with the open book $(S_{g,r},\phi)$ is tight if and only if the fractional Dehn twists coefficient of $\phi$ around every boundary component is non-negative \cite{ch}. In this case, the contact structure is also Stein fillable \cite{ch}. So, in particular, if $(M,\xi)$ is supported by an open book with periodic monodromy, then $\xi$ is tight if and only if $c^+(\xi) \neq 0$. Therefore, Theorem \ref{thm:periodic} may be reformulated as follows.

\begin{theorem}
\label{thm:periodic2}
Suppose that $(M,\xi)$ is supported by a genus one open book with $r$ binding components and periodic monodromy. If $c^+(\xi) \neq 0$, then $r\geq 1+4gr(c^+(\xi))$. 
\end{theorem}

To prove this theorem, we bound the grading shifts associated to the maps induced by capping off, and we use the fact that $gr(c^+(\xi)) \leq 0$ whenever $\xi$ is tight and is supported by a genus one open book with \emph{one} boundary component and periodic monodromy (see Table \ref{table:grfdtc}). Before we compute these grading shifts, we must know that they are well-defined. To this end, we establish the following.

\begin{proposition}
\label{prop:tor}
Suppose the contact 3-manifold $(M,\xi)$ is supported by an open book $(S,\phi)$ for which $\phi^m$ is freely isotopic to the identity. If the fractional Dehn twist coefficients of $\phi$ are non-negative, then $c_1(\spt_{\xi})$ is a torsion class.
\end{proposition}

\begin{proof}[Proof of Proposition \ref{prop:tor}]
It suffices to show that $\langle c_1(\spt_{\xi}), h \rangle = 0$ for every $h\in H_2(M,\zz)$. Let $(\Sigma, \alpha, \beta, z)$ be a standard pointed Heegaard diagram for the open book $(S_{g,r},\phi)$, and let $\mathcal{D}_1, \dots, \mathcal{D}_k$ denote the connected regions of $\Sigma \setminus \cup \alpha_i \setminus \cup \beta_i$. Recall that a doubly-periodic domain for this pointed Heegaard diagram is a 2-chain $\mathcal{P} = \sum_i a_i \mathcal{D}_i$ whose boundary is a sum of $\alpha$ and $\beta$ curves, and whose multiplicity at the basepoint $z$ is 0. It is often convenient to think of a periodic domain as a linear relation in $H_1(\Sigma;\zz)$ amongst the $\alpha$ and $\beta$ curves. Doubly-periodic domains are in one-to-one correspondence with elements of $H_2(M;\zz)$; we denote by $H(\mathcal{P})$ the homology element corresponding to $\mathcal{P}$. Suppose that $\mathbf{y} = \{y_1,\dots,y_n\}$ is the intersection point between $\Tb$ and $\Ta$ described in Subsection \ref{ssec:heeg} for which $[\mathbf{y},0]\in \cfp(\Sigma,\beta,\alpha,z)$ represents $c^+(\xi)$ (here, $n=2g+r-1$). Then $\spt_{\xi}$ is the $Spin^c$ structure associated to $\mathbf{y}$; that is, $\spt_{\xi} = \spc_z(\mathbf{y})$ \cite{osz1}. So, our goal is to show that $\langle c_1(\spc_z(\mathbf{y})), H(\mathcal{P}) \rangle = 0$ for every doubly-periodic domain $\mathcal{P}$.

The \emph{Euler measure} of a region $\mathcal{D}_i$ is the quantity $$\widehat{\chi}(\mathcal{D}_i) = \chi(\text{int }\mathcal{D}_i) - \frac{1}{4}(\# \text{corner points of }\mathcal{D}_i),$$ where corner points are to be counted with multiplicity \cite{osz5}. We extend the definition of Euler measure to 2-chains linearly. Let $n_{\mathbf{y}}(\mathcal{P})$ be the sum of the local multiplicities of $\mathcal{P}$ at the points $y_i \in \mathbf{y}$. By \cite[Proposition 7.5]{osz14}, $$\langle c_1(\spc_z(\mathbf{y})), H(\mathcal{P}) \rangle = \widehat{\chi}(\mathcal{P}) + 2n_{\mathbf{y}}(\mathcal{P}).$$ Below, we prove that both $\widehat{\chi}(\mathcal{P})$ and $n_{\mathbf{y}}(\mathcal{P})$ vanish for every doubly-periodic domain $\mathcal{P}$.  

Suppose that the arcs $a_1,\dots,a_n$ on $S=S_{g,r}$, used to form the $\alpha$ and $\beta$ curves, are those shown in Figure \ref{fig:page}. Let $B_i$ be the boundary component of $S$ which intersects the arcs $a_i$ and $a_{i-1}$ (unless $i=1$, in which case $B_1$ is the boundary component which intersects only $a_1$). For $i=1,\dots,r$, let $d_i$ be the oriented curve on $\Sigma = S_{1/2} \cup -S_0$ defined by $d_i = B_i \times \{1/2\},$ where $d_i$ inherits its orientation from the boundary orientation on $B_i$. We orient the $\alpha$ and $\beta$ curves so that the orientation of the arc $\alpha_i \cap S_{1/2}$ agrees with that of $\beta_i \cap S_{1/2}$. Furthermore, we require that $\alpha_i \cdot d_i = +1$ for $i=1,\dots,r-1$. 

We may assume that for some fixed integer $m$, the fractional Dehn twist coefficient of $\phi$ around each $B_i$ is given by $k_i/m$ for some integer $k_i\geq 0$. Then, $\phi^m$ is isotopic to a product of Dehn twists $t_{c_1}^{k_1}\cdots t_{c_r}^{k_r}$, as discussed at the beginning of this section. Recall that the arc $b_i$ on $S$ is obtained from $a_i$ via a small isotopy, as described in Subsection \ref{ssec:heeg}. Let $b^{(1)}_i$ denote the arc $b_i$, and let $b_i^{(j)}$ be the arc on $S$ obtained from $b_i^{(j-1)}$ via a similar isotopy for each $j=2,\dots,m$ (so that $b_i^{(j)}$ intersects $b_i^{(j-1)}$ transversely in one point and with positive sign). Recall that $\alpha_i$ and $\beta_i$ are defined by $$\alpha_i = a_i \times \{1/2\} \cup a_i \times \{0\},$$ $$\beta_i = b_i \times \{1/2\} \cup \phi(b_i) \times \{0\}.$$ For $j=2,\dots,m$, we define $$\beta^{(j)}_i = b^{(j)}_i \times \{1/2\} \cup \phi^j(b^{(j)}_i) \times \{0\}.$$

Suppose $\mathcal{P}$ is a doubly-periodic domain specified by the relation \begin{equation}\label{eqn:rel1}\sum_i s_i\alpha_i + \sum_i t_i \beta_i = 0\end{equation} in $H_1(\Sigma;\zz)$. Let $\delta_i$ be a curve on $S$ which intersects the arc $a_i$ exactly once (and does not intersect the other $a_j$). The curve $\delta_i \times \{1/2\} \subset \Sigma$ must algebraically intersect $\partial \mathcal{P}$ zero times; that is, $(\delta_i \times \{1/2\})\cdot \partial \mathcal{P} = \pm (s_i+t_i)=0$.  We can therefore express the relation in Equation \ref{eqn:rel1} by \begin{equation}\label{eqn:per}\sum_i s_i (\alpha_i-\beta_i)=0.\end{equation} But this implies that $n_{\mathbf{y}}(\mathcal{P})=0$ (see Figure \ref{fig:local} for the local picture of $\mathcal{P}$ near the intersection point $y_i$). To see that $\widehat{\chi}(\mathcal{P}) = 0$ as well, we consider the pointed Heegaard multi-diagram $(\Sigma, \alpha,\beta,\beta^{(2)}, \dots, \beta^{(m)},z)$. 

\begin{figure}[!htbp]

\labellist 
\hair 2pt 
\small\pinlabel $\bullet$ at 94 88
\pinlabel $\alpha_i$ at 125 175
\pinlabel $\beta_i$ at 64 176
\pinlabel $\bullet z$ at 151 30

\pinlabel $S_{1/2}$ at 175 139
\pinlabel $-S_{0}$ at 175 171

\pinlabel $\pm s_i$ at 94 126
\pinlabel $\mp s_i$ at 94 50

\pinlabel $0$ at 130 88
\pinlabel $0$ at 54 88
\tiny\pinlabel $y_i$ at 107 89

\endlabellist 

\begin{center}
\includegraphics[width=3.9cm]{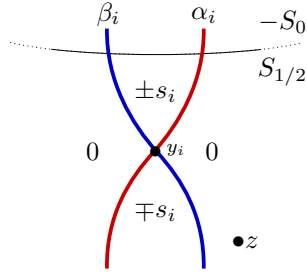}
\caption{\quad The coefficients of $\mathcal{P}$ near the intersection point $y_i$.}
\label{fig:local}
\end{center}
\end{figure}

The relation in Equation \ref{eqn:per} implies that \begin{equation}\label{eqn:per2}\sum_i s_i (\beta^{(j-1)}_i-\beta^{(j)}_i)=0\end{equation} in $H_1(\Sigma;\zz)$ as well, for each $j=2,\dots,m$. Let $\mathcal{P}_j$ be the doubly-periodic domain specified by the relation in Equation \ref{eqn:per2}. The doubly-periodic domain $\mathcal{P}_{sum} = \mathcal{P} + \mathcal{P}_2 + \dots + \mathcal{P}_m$ is therefore specified by the relation \begin{equation}\label{eqn:per3}\sum_i s_i(\alpha_i - \beta^{(m)}_i)=0\end{equation} obtained by summing the relation in Equation \ref{eqn:per} with those in Equation \ref{eqn:per2}. Since $\widehat{\chi}(\mathcal{P}_j)= \widehat{\chi}(\mathcal{P})$ for each $j=2,\dots,m$, and Euler measure is additive, \begin{equation}\label{eqn:Psum}\widehat{\chi}(\mathcal{P}_{sum}) = m\widehat{\chi}(\mathcal{P}).\end{equation} 

Observe that $(\Sigma, \alpha, \beta^{(m)},z)$ is a standard pointed Heegaard diagram for $(S_{g,r},\phi^m)$, and recall that $\phi^m$ is isotopic to $t_{c_1}^{k_1}\cdots t_{c_r}^{k_r}$. Then, in $H_1(\Sigma;\zz)$, \begin{equation}\label{eqn:poss}\alpha_i - \beta^{(m)}_i = 
\left\{\begin{array}{lll} -k_id_i+k_{i+1}d_{i+1},& 1 \leq i < r,\\ 
0,  & r \leq i \leq n.
 \end{array} \right.\end{equation} So, if the relation in Equation \ref{eqn:per3} holds, then $\sum_{i<r} s_i(-k_id_i+k_{i+1}d_{i+1})=0$ in $H_1(\Sigma;\zz)$ as well. But any relation in $H_1(\Sigma;\zz)$ amongst the curves $d_1,\dots,d_r$ is of the form $t(d_1+\dots + d_r)=0$. Hence, $$\sum_{i<r} s_i(-k_id_i+k_{i+1}d_{i+1}) = t(d_1+\dots + d_r).$$ On the other hand, since all of the $k_i$ are non-negative, this can only happen if $t=0$. It follows that $s_i = s_{i-1}$ if $k_i\neq 0$ (unless $i=1$, in which case $k_1\neq 0$ implies that $s_1=0$).
 
Therefore, the relation in Equation \ref{eqn:per3} breaks up into smaller relations of the form \begin{equation}\label{eqn:per4}s_i(\alpha_i-\beta^{(m)}_i)=0,\end{equation} for $i\geq r$, and \begin{equation}\label{eqn:per5}\sum_{i_1\leq i \leq i_2<r}s_i(\alpha_i - \beta^{(m)}_i)= 0,\end{equation} where $s_i = s_j$ for $i$ and $j$ between $i_1$ and $i_2$. It is not hard to see directly that the doubly-periodic domains given by the relations in Equations \ref{eqn:per4} and \ref{eqn:per5} have vanishing Euler measure. In either case, these periodic domains, thought of as linear combinations of regions in $\Sigma \setminus \cup \alpha_i \setminus \cup \beta^{(m)}_i$, each consist of two canceling bigon regions together with square regions (whose Euler measures are zero). See Figures \ref{fig:perdom2} and \ref{fig:perdom} for reference. It follows that $\widehat{\chi}(\mathcal{P}_{sum}) = 0$, which, in turn, implies that $\widehat{\chi}(\mathcal{P}) = 0$, by Equation \ref{eqn:Psum}. This completes the proof of Proposition \ref{prop:tor}.

\begin{figure}[!htbp]

\labellist 
\hair 2pt 
\small\pinlabel $\dots$ at -28 184
\pinlabel $-S_0$ at -28 234
\pinlabel $S_{1/2}$ at -28 134

\pinlabel $\alpha_i$ at 275 145
\pinlabel $\beta^{(m)}_i$ at 305 85

\pinlabel $\bullet z$ at 55 60

\endlabellist 

\begin{center}
\includegraphics[width=4.8cm]{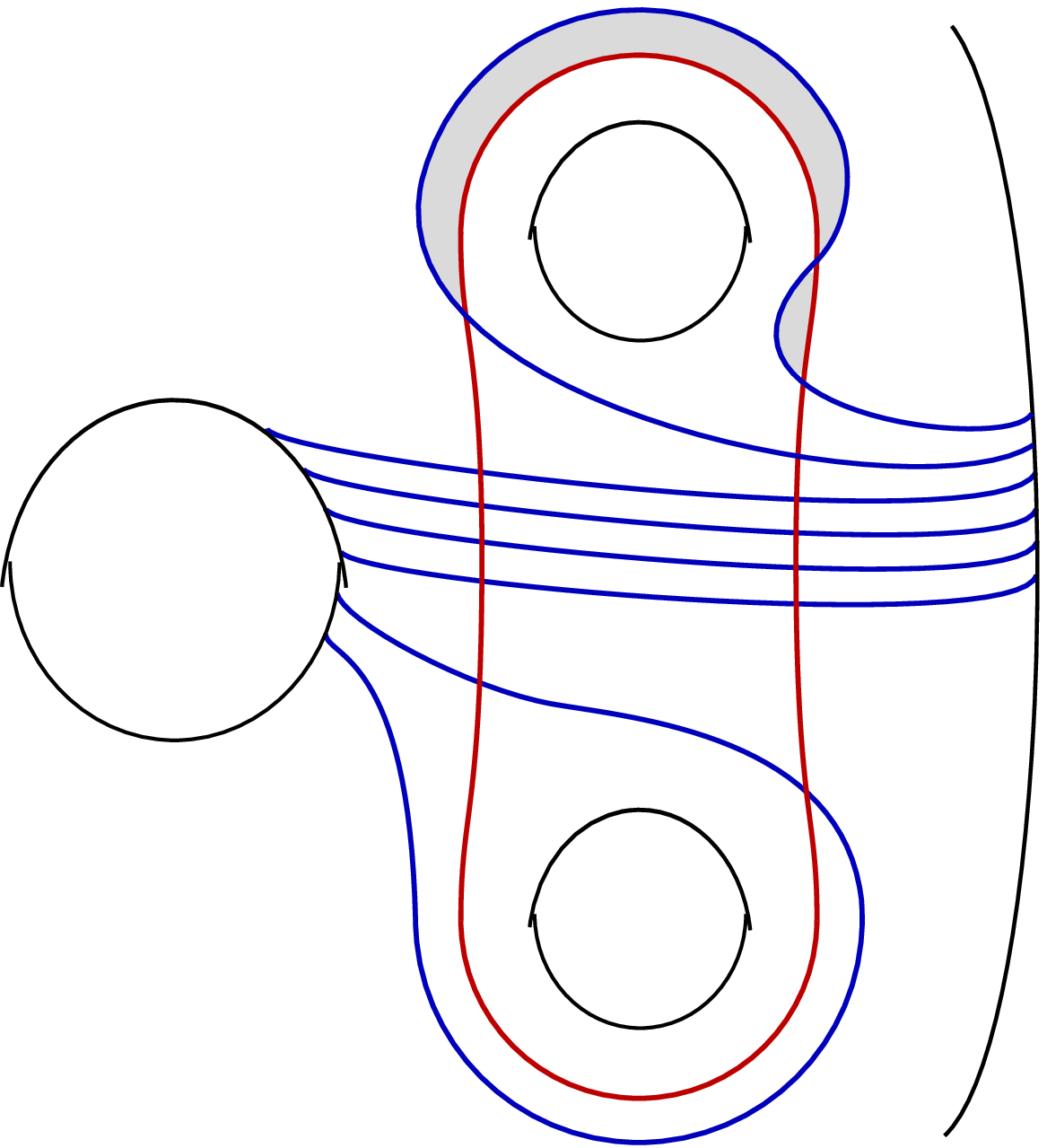}
\caption{\quad Shown here are $\alpha_i$ and $\beta^{(m)}_i$ for some $i\geq r$. The region bounded by these curves is a periodic domain corresponding to a relation as in Equation \ref{eqn:per4}. Note that it consists of square regions and two canceling bigon regions (which we have shaded). In this figure, $k_r=3$.}
\label{fig:perdom2}
\end{center}
\end{figure}

\begin{figure}[!htbp]

\labellist 
\hair 2pt 
\small\pinlabel $\bullet z$ at 150 28
\pinlabel $\dots$ at 560 87
\pinlabel $-S_0$ at 560 125
\pinlabel $S_{1/2}$ at 560 49

\pinlabel $\beta^{(m)}_{i_1}$ at 67 0
\pinlabel $\beta^{(m)}_{i_1+1}$ at 258 0
\pinlabel $\beta^{(m)}_{i_1+2}$ at 440 0
\pinlabel $\beta^{(m)}_{i_2}$ at 687 0
\endlabellist 

\begin{center}
\includegraphics[width=12.9cm]{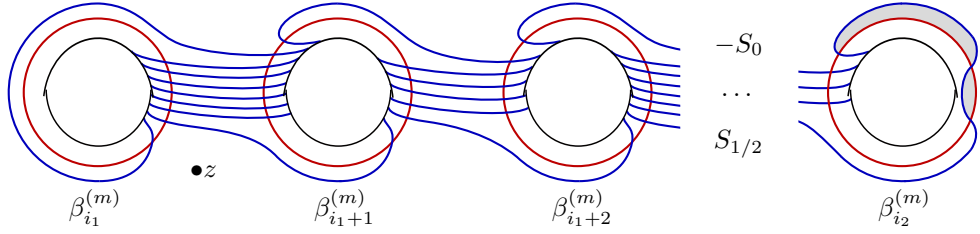}
\caption{\quad The region bounded by these $\alpha$ and $\beta^{(m)}$ curves is a periodic domain corresponding to a relation as in Equation \ref{eqn:per5}. It consists of square regions and two canceling bigon regions (which we have shaded). In this figure, $k_{i_1} = 0$, $k_{i_1+1} = 4$, $k_{i_1+2} = 3$, $k_{i_1+3} = 4$, $k_{i_2} = 2$ and $k_{i_2+1}=0$. }
\label{fig:perdom}
\end{center}
\end{figure}
\end{proof}

\subsection{Grading shifts and the proof of Theorem \ref{thm:periodic}} Below, we study the grading shifts associated to the maps induced by capping off. Suppose $(S_{g,r-1},\phi')$ is the open book obtained from $(S_{g,r},\phi)$ by capping off one of the boundary components of $S_{g,r}$. Let $W$ be the corresponding 2-handle cobordism from $-M_{S_{g,r-1},\phi'}$ to $-M_{S_{g,r},\phi}$. If $\phi$ is periodic with non-negative fractional Dehn twist coefficients, then the same is true of $\phi'$, and it follows from Proposition \ref{prop:tor} that the contact invariants $c^+(S_{g,r-1},\phi')$ and $c^+(S_{g,r},\phi)$ have well-defined $\mathbb{Q}$-gradings. Since $F^+_{W,\spc_0}$ sends $c^+(S_{g,r-1},\phi')$ to $c^+(S_{g,r},\phi),$ by Theorem \ref{thm:nat}, the grading shift formula in \cite{osz5} gives \begin{equation}\label{eqn:grshift}gr(c^+(S_{g,r},\phi)) - gr(c^+(S_{g,r-1},\phi')) =\frac{c_1(\spc_0)^2-2\chi(W)-3\sigma(W)}{4}.\end{equation} 

\begin{lemma}
\label{lem:negdef}
The cobordism $W$ either has trivial intersection form, or $b_2^+(W)=0$ and $\sigma(W) = -1$.
\end{lemma}

Since $W$ is obtained from a single 2-handle attachment, $\chi(W)=1$. Together with Lemma \ref{lem:negdef}, this implies that \begin{equation}\label{eqn:gr}gr(c^+(S_{g,r},\phi)) - gr(c^+(S_{g,r-1},\phi')) \leq 1/4.\end{equation} Suppose that after capping off all but one of the boundary components of $S_{g,r}$, we are left with an open book $(S_{g,1},\phi'')$. It follows from Inequality \ref{eqn:gr} that \begin{equation*}\label{eqn:gr2}gr(c^+(S_{g,r},\phi)) \leq (r-1)/4 +gr(c^+(S_{g,1},\phi'')).\end{equation*} Now, suppose $g=1$. As mentioned in the beginning of this section, $gr(c^+(S_{1,1},\phi'')) \leq 0$ (see Table \ref{table:grfdtc}). As a result, $$gr(c^+(S_{1,r},\phi)) \leq (r-1)/4,$$ which is equivalent to the statement of Theorem \ref{thm:periodic2}. All that remains is to prove Lemma \ref{lem:negdef}.

\begin{proof}[Proof of Lemma \ref{lem:negdef}]
Let $\Delta$ denote the 2-simplex with edges, $e_{\beta}$, $e_{\gamma}$ and $e_{\alpha},$ described in Subsection \ref{ssec:cap}. The pointed Heegaard triple-diagram $(\Sigma,\beta,\gamma,\alpha,z)$ associated to the capping off operation above (and defined in Subsection \ref{ssec:heeg}) specifies an identification space $$X_{\beta,\gamma,\alpha} = \frac{(\Delta \times \Sigma) \amalg (e_{\beta} \times H_{\beta}) \amalg (e_{\gamma} \times H_{\gamma}) \amalg (e_{\alpha} \times H_{\alpha})}{(e_{\beta}\times \Sigma) \sim (e_{\beta} \times \partial H_{\beta}),\, (e_{\gamma}\times \Sigma) \sim (e_{\gamma} \times \partial H_{\gamma}),\, (e_{\alpha}\times \Sigma) \sim (e_{\alpha} \times \partial H_{\alpha})},$$ 
where $H_{\beta}$, $H_{\gamma}$ and $H_{\alpha}$ are the handlebodies corresponding to the $\beta$, $\gamma$ and $\alpha$ curves (see the diagram on the left of Figure \ref{fig:Xwhitney} for a schematic picture of $X_{\beta,\gamma,\alpha}$). After rounding corners, $X_{\beta,\gamma,\alpha}$ is a smooth cobordism with boundary $-M_{\beta,\gamma} - M_{\gamma, \alpha} + M_{\beta,\alpha}$. (Here, $M_{\beta,\gamma}$ is the 3-manifold specified by the Heegaard diagram $(\Sigma,\beta,\alpha)$, and similarly for $M_{\gamma, \alpha}$ and $M_{\beta,\alpha}$.) In fact, $X_{\beta,\gamma,\alpha}$ is just the complement of a neighborhood of a 1-complex in the cobordism $W$, so the intersection form of $X_{\beta,\gamma,\alpha}$ is the same as that of $W$ (refer to \cite{osz8,osz5} for more details).

Elements of $H_2(X_{\beta,\gamma,\alpha};\zz)$ are in one-to-one correspondence with triply-periodic domains for the Heegaard diagram $(\Sigma,\beta,\gamma,\alpha,z)$. If $\mathcal{P}$ is a triply-periodic domain, we denote the $\beta$, $\gamma$ and $\alpha$ components of $\partial \mathcal{P}$ by $\partial_{\beta} \mathcal{P}$, $\partial_{\gamma} \mathcal{P}$ and $\partial_{\alpha} \mathcal{P}$. The homology class corresponding to a triply-periodic domain $\mathcal{P}$ is constructed as follows. Pick a point $p\in \Delta$, and consider the copy of $\mathcal{P}$ contained in $\{p\} \times \Sigma.$ Attach cylinders, connecting each component of $\partial_{\beta} \mathcal{P}$ in $\{p\} \times \Sigma$ with the corresponding component in $\{u\}\times \Sigma$ for some $u \in e_{\beta}$. Then cap off these boundary components with disks inside $\{u\} \times H_{\beta}$. Do the same for the components of $\partial_{\gamma} \mathcal{P}$ and $\partial_{\alpha} \mathcal{P}$. We denote this homology class by $H(\mathcal{P})$. 

The middle diagram in Figure \ref{fig:Xwhitney} shows a schematic picture of this construction. The point labeled $p$ represents the copy of $\mathcal{P}$ in $\{p\} \times \Sigma$, and the three legs represent the attaching cylinders for the components of $\partial_{\beta} \mathcal{P}$, $\partial_{\gamma} \mathcal{P}$ and $\partial_{\alpha} \mathcal{P}$. The rightmost diagram is meant to represent the intersection of two such homology classes, $H(\mathcal{P})$ and $H(\mathcal{P}')$. The $\beta$ attaching cylinders of $H(\mathcal{P})$ intersect the $\gamma$ attaching cylinders of $H(\mathcal{P}')$ at points in $\{q\}\times \Sigma$, and it's not hard to check that the algebraic intersection number \begin{equation}\label{eqn:selfint}H(\mathcal{P}) \cdot H(\mathcal{P'})= (\partial_{\beta}\mathcal{P}) \cdot (\partial_{\gamma}\mathcal{P}').\end{equation} In particular, note that if $\mathcal{P}$ is a doubly-periodic domain (by which we mean that $\partial \mathcal{P}$ consists of only two of the three types of attaching curves) then $H(\mathcal{P})$ pairs trivially with every element in $H_2(X_{\beta,\gamma,\alpha};\zz)$.

\begin{figure}[!htbp]

\labellist 
\hair 2pt 
\small
\pinlabel $\Delta\times \Sigma$ at 129 96
\pinlabel  \rotatebox{59}{$e_{\beta}\times H_{\beta}$} at 65 128
\pinlabel \rotatebox{-60}{$e_{\gamma}\times H_{\gamma}$} at 194 128
\pinlabel $e_{\alpha}\times H_{\alpha}$ at 129 25

\pinlabel $p$ at 425 68
\pinlabel $p$ at 731 68
\pinlabel $p'$ at 799 75

\pinlabel $q$ at 767 105

\endlabellist 

\begin{center}
\includegraphics[width=13cm]{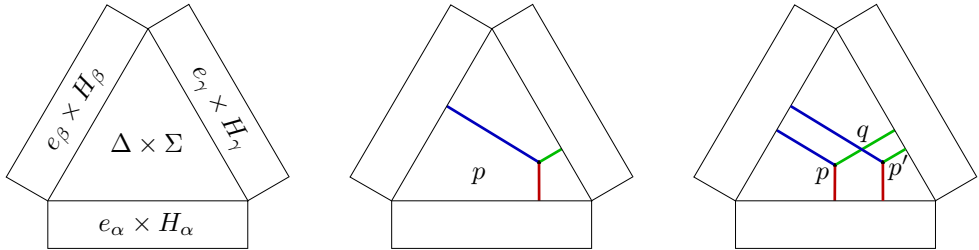}
\caption{\quad From left to right: the cobordism $X_{\beta,\gamma,\alpha}$, a homology class $H(\mathcal{P})$, and the intersection of two classes, $H(\mathcal{P})$ and $H(\mathcal{P'})$.}
\label{fig:Xwhitney}
\end{center}
\end{figure}

Suppose that the intersection form of $X_{\beta,\gamma,\alpha}$ is non-trivial. Let $\mathcal{P}$ be a triply-periodic domain given by the relation $$N\gamma_1 + \sum_i s_i \alpha_i + \sum_i t_i \beta_i = 0$$ in $H_1(\Sigma;\zz)$, where $N\neq 0$. Since $\gamma_i$ is isotopic to $\beta_i$ for $i=2,\dots,n$, every triply-periodic domain differs from some multiple of $\mathcal{P}$ by a sum of doubly-periodic domains. Let $d_1,\dots, d_r$ be the curves defined in the proof of Proposition \ref{prop:tor}, and orient the $\alpha$ and $\beta$ curves as before. We orient $\gamma_1$ in the same direction as $d_1$. The same argument used in the proof of Proposition \ref{prop:tor} shows that $r_i = -s_i$ for each $i$, so $\mathcal{P}$ is given by the relation \begin{equation}\label{eqn:rela}N\gamma_1 + \sum_i s_i (\alpha_i -  \beta_i) = 0.\end{equation} By Equation \ref{eqn:selfint}, $$H(\mathcal{P})^2= (\partial_{\beta}\mathcal{P}) \cdot (\partial_{\gamma}\mathcal{P}) = -Ns_1(\beta_1 \cdot \gamma_1) = -Ns_1,$$ so Lemma \ref{lem:negdef} follows if we can show that $Ns_1>0$. 

As before, we assume that the fractional Dehn twist coefficient of $\phi$ around $B_i$ is $k_i/m$, where $k_i\geq 0$. The relation in Equation \ref{eqn:rela} implies that \begin{equation}\label{eqn:rela2}N\gamma_1 + \sum_i s_i (\beta^{(j-1}_i -  \beta^{(j)}_i) = 0\end{equation} in $H_1(\Sigma;\zz)$ for $j=2,\dots,m$. Adding the relations in Equations \ref{eqn:rela} and \ref{eqn:rela2}, we find that \begin{equation}\label{eqn:rela3}mN\gamma_1 + \sum_i s_i (\alpha_i -  \beta^{(m)}_i) = 0\end{equation} in  $H_1(\Sigma;\zz)$. After making the substitutions from Equation \ref{eqn:poss}, and noting that $\gamma_1 = d_1$, it follows that $$mNd_1 + \sum_{i<r} s_i(-k_id_i+k_{i+1}d_{i+1})=0$$ in $H_1(\Sigma;\zz)$ as well. As in the proof of Proposition \ref{prop:tor}, this implies that \begin{equation}\label{eqn:relt}mNd_1 + \sum_{i<r} s_i(-k_id_i+k_{i+1}d_{i+1}) = t(d_1+\dots + d_r)\end{equation} for some $t$. If $t=0$, then $mN=s_1k_1.$ We are assuming that $X_{\beta,\gamma,\alpha}$ has non-trivial intersection form, so $H(\mathcal{P})^2= - Ns_1 \neq 0$. Therefore, $k_1$ is strictly greater than zero, and it follows that $Ns_1>0$, as hoped. If $t\neq 0$, we can assume, without loss of generality, that $t>0$. Then all $k_i$ and $s_i$ are strictly greater than zero, and $mN-s_1k_1>0$, which implies that $Ns_1>0$ as well.
\end{proof}

\subsection{A formula for the $d_3$ invariant}
Below, we explicitly compute the grading shift in Equation \ref{eqn:grshift} in terms of the fractional Dehn twist coefficients of $\phi$. If the intersection form of $W$ is trivial, then this grading shift is simply $-1/2$. Otherwise, $b_2^+(W)=0$ and $\sigma(W)=-1$ (by Lemma \ref{lem:negdef}), and the grading shift is $(c_1(\spc_0)^2 +1)/4$. According to Proposition \ref{prop:tor}, $c_1(\spc_0)$ is sent to zero by the restriction map $$H^2(W;\mathbb{Q})\rightarrow H^2(\partial W;\mathbb{Q}).$$ Therefore, $c_1(\spc_0)$ is the image of a class $k \cdot PD(\lambda)$ under the map $$H^2(W,\partial W;\mathbb{Q}) \rightarrow H^2(W;\mathbb{Q}),$$ where $\lambda$ is a generator of the dimension one subspace of elements $B_2^-(W)\subset H_2(W)$ with negative self-intersection. By definition, $$c_1(\spc_0)^2 = (k \cdot PD(\lambda))^2 = k^2 \cdot \lambda^2 = \frac{\langle c_1(\spc_0), \lambda \rangle ^2}{\lambda^2}.$$ Recall that $\spc_0 = \spc_z(\psi_0)$, where $\psi_0$ is the homotopy class of Whitney triangles whose domain $\mathcal{D}(\psi_0)$ is $\Delta_1 + \dots + \Delta_n$, and let $H(\mathcal{P})$ be a class which generates $B_2^-(X_{\beta,\gamma,\alpha})$. Then, the equation above becomes \begin{equation}\label{eqn:csquared}c_1(\spc_0)^2 = \frac{\langle c_1(\spc_z(\psi_0)), H(\mathcal{P}) \rangle ^2}{H(\mathcal{P})^2}.\end{equation}

To compute $\langle c_1(\spc_z(\psi_0)), H(\mathcal{P}) \rangle$, we introduce the \emph{dual spider number} of a Whitney triangle $$u:\Delta \rightarrow Sym^n(\Sigma)$$ and a triply-periodic domain $\mathcal{P}$, following the exposition in \cite{osz5}. First, note that the orientations on the $\beta$, $\gamma$ and $\alpha$ curves induce ``leftward" pointing normal vector fields along the curves. Let $\beta_i'$, $\gamma_i'$ and $\alpha_i'$ be copies of the attaching curves $\beta_i$, $\gamma_i$ and $\alpha_i$, translated slightly in these normal directions, and let $\mathbb{T}_{\beta'}$, $\mathbb{T}_{\gamma'}$ and $\mathbb{T}_{\alpha'}$ denote the corresponding tori in $Sym^n(\Sigma)$. By construction, $u(e_{\beta})$ misses $\mathbb{T}_{\beta'}$, $u(e_{\gamma})$ misses $\mathbb{T}_{\gamma'}$ and $u(e_{\alpha})$ misses $\mathbb{T}_{\alpha'}$.

Let $x$ be an interior point of $\Delta$ so that $u(x)$ misses the $\beta'$, $\gamma'$ and $\alpha'$ curves, and choose three oriented paths, $b$, $c$ and $a$, from $x$ to $e_{\beta}$, $e_{\gamma}$ and $e_{\alpha}$, respectively. Let $\partial_{\beta'}(\mathcal{P})$, $\partial_{\gamma'}(\mathcal{P})$ and $\partial_{\alpha'}(\mathcal{P})$ be the 1-chains obtained by translating the corresponding components of $\partial \mathcal{P}$ in the normal directions described above. The dual spider number of $u$ and $\mathcal{P}$ is given by $$\sigma(u,\mathcal{P}) = n_{u(x)}(\mathcal{P}) + \partial_{\beta'}(\mathcal{P}) \cdot b + \partial_{\gamma'}(\mathcal{P}) \cdot c + \partial_{\alpha'}(\mathcal{P}) \cdot a.$$ In \cite{osz5}, Ozsv{\'a}th and Szab{\'o} prove that $$\langle c_1(\spc_z(\psi_0)), H(\mathcal{P})\rangle = \widehat{\chi}(\mathcal{P}) + \#(\partial \mathcal{P}) + 2\sigma(u,\mathcal{P})$$ for any Whitney triangle $u$ representing $\mathcal{P}$.

Suppose $\phi$ is a periodic diffeomorphism of $S=S_{g,r}$ with fractional Dehn twist coefficients $0\leq k_1/m \leq \dots \leq k_r/m,$ and suppose that the intersection form of $W$ is non-trivial. Let $\mathcal{P}$ be the triply-periodic domain specified by the relation in Equation \ref{eqn:rela}. To compute $\widehat{\chi}(\mathcal{P})$, we consider the pointed Heegaard multi-diagram $(\Sigma,\alpha,\beta,\beta^{(2)},\dots,\beta^{(m)},\gamma,z)$. Let $\mathcal{P}_j$ be the triply-periodic domain specified by the relation in Equation \ref{eqn:rela2} for $j=2,\dots,m$. Then $\mathcal{P}_{sum} = \mathcal{P} + \mathcal{P}_2 + \dots + \mathcal{P}_m$ is specified by the relation in Equation \ref{eqn:rela3}, and $$\widehat{\chi}(\mathcal{P}_{sum}) = m \widehat{\chi}(\mathcal{P}),$$ as before. Per Equation \ref{eqn:poss}, this relation breaks up into relations in $H_1(\Sigma;\mathbb{Z})$ of the form $$s_i(\alpha_i-\beta_i^{(m)}) = 0,$$ for $i\geq r,$ and \begin{equation}\label{eqn:relas}mNd_1 + \sum_{i<r} s_i(-k_id_i+k_{i+1}d_{i+1})=0.\end{equation} As noted previously, the doubly-periodic domains specified by the former relations have Euler measure zero, and the latter relation implies that $$mNd_1 + \sum_{i<r} s_i(-k_id_i+k_{i+1}d_{i+1}) = t(d_1+\dots + d_r)$$ for some $t$. Suppose that $t=0$. If no $k_i$ is zero, then all of the $s_i$ must vanish. But this implies that $H(\mathcal{P})^2 = -Ns_1 = 0$, which contradicts our assumption on the intersection form of $W$. If some $k_i$ vanishes, then $k_1$ must vanish since $0\leq k_1 \leq k_i$ by assumption. But this too implies that $H(\mathcal{P})^2 = 0.$ So, it must be the case $t\neq 0$ and $k_i>0$ for all $i$. We may assume, without loss of generality, that $t=-k_r$. Then \begin{equation}\label{eqn:s}s_i = -k_r(1/k_r + 1/k_{r-1} + \dots + 1/k_{i+1}). \end{equation} We define $s_0$ using this formula as well; note that $s_0 = mN/k_1$.

\begin{figure}[!htbp]

\labellist 
\hair 2pt 
\small\pinlabel $\bullet z$ at 240 45
\pinlabel $\dots$ at 460 127
\pinlabel $-S_0$ at 460 162
\pinlabel $S_{1/2}$ at 460 89

\pinlabel $F$ at 240 230
\pinlabel $\beta^{(m)}_{1}$ at 150 37
\pinlabel $\beta^{(m)}_{2}$ at 333 37
\pinlabel $\beta^{(m)}_{r-1}$ at 600 37

\tiny\pinlabel $\dots$ at 787 229
\pinlabel $\dots$ at 787 27

\endlabellist 

\begin{center}
\includegraphics[width=14.2cm]{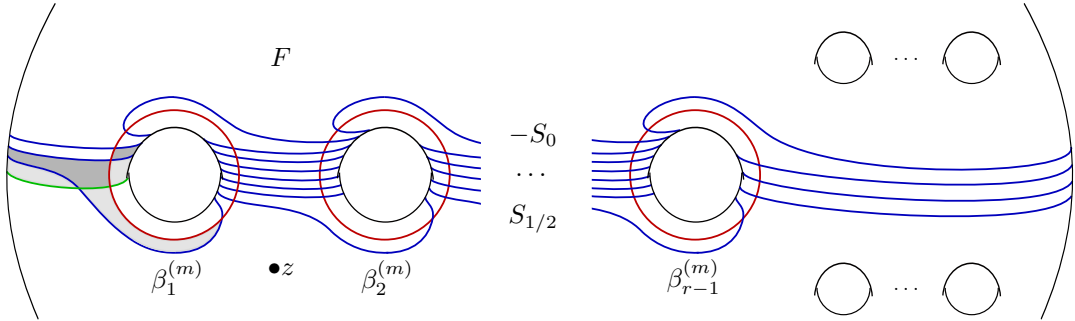}
\caption{\quad A portion of the surface $\Sigma$. The region bounded by the $\alpha$, $\beta^{(m)}$ and $\gamma$ curves is the periodic domain $\mathcal{P}_{sum}'$. It is composed of squares, two triangular regions (gray), a pentagonal region (dark gray) and the region $F$.}
\label{fig:perdom3}
\end{center}
\end{figure}

The triply-periodic domain $\mathcal{P}_{sum}'$, given by the relation in Equation \ref{eqn:relas}, is composed of square regions, two triangular regions, a pentagonal region and a region $F$ which has genus $g$, one boundary component, and $4(r-1)$ corners (see Figure \ref{fig:perdom3}). It is easy to check that the contributions of the triangular regions and the pentagonal region to $\widehat{\chi}(\mathcal{P}_{sum'})$ cancel. Since the region $F$ has multiplicity $k_r$ in $\mathcal{P}_{sum}'$, $$\widehat{\chi}(\mathcal{P}_{sum})= \widehat{\chi}(\mathcal{P}_{sum}') = k_r(2-2g-r),$$ and, hence, \begin{equation}\label{eqn:emeasure}\widehat{\chi}(\mathcal{P}) = k_r(2-2g-r)/m.\end{equation}

Let $u:\Delta \rightarrow Sym^n(\Sigma)$ be a representative of the homotopy class $\psi_0$. The local contribution of $\Delta_i$ to the dual spider number $\sigma(u,\mathcal{P})$ is $-|s_i|$. On the other hand, the number of boundary components of $\mathcal{P}$ is $|N| + |s_1| + \dots + |s_n|$. So, the quantity $$\langle c_1(\spc_z(\psi_0)), H(\mathcal{P})\rangle = \widehat{\chi}(\mathcal{P}) + \#(\partial \mathcal{P}) + 2\sigma(u,\mathcal{P})$$ is simply $|N| = -N = -k_1s_0/m$. And we saw in the previous subsection that $H(\mathcal{P})^2 = -Ns_1=k_1s_0s_1/m$. As a result, \begin{equation}\label{eqn:grshift2}c_1(\spc_0)^2 = \frac{(k_r(2-2g-r)/m - k_1s_0/m)^2}{k_1s_0s_1/m} = \frac{(k_r(2-2g-r) - k_1s_0)^2}{mk_1s_0s_1},\end{equation} by Equation \ref{eqn:csquared}. 

We now restrict our focus to genus one open books. Suppose that $\xi$ is a tight contact structure supported by an open book $(S_{1,1},\phi)$ with periodic $\phi$. The table below lists the grading of $c^+(\xi)$ as a function of the fractional Dehn twist coefficient (FDTC) of $\phi$. This follows from the grading calculations in \cite[Section 6]{bald4} (for non-integral FDTC's) and \cite[Proposition 9]{jabmark2} (for integral FDTC's).

\begin{table}[ht]
 \label{table:grfdtc}
\begin{center}
\begin{tabular}{|c|c|c|} \hline 
$FDTC$ & $gr(c^+(\xi))$ \cr \hline \hline 
$(6k+5)/6$  &$-2$   \cr \hline 
$(4k+3)/4$  & $-7/4$  \cr \hline 
$(3k+2)/3$  & $-3/2$  \cr \hline
$(6k+1)/6$  & $-1/2$   \cr \hline
$(4k+1)/4$  & $-1/4$  \cr \hline
$(3k+1)/3$  & $0$   \cr \hline
$(2k-1)/2$  & $-1$   \cr \hline
$k$  & $-1$   \cr \hline
\end{tabular}
\end{center}
\vspace{3mm}
 \caption{\quad Grading versus fractional Dehn twist coefficient. In this table, $k\geq 0$; otherwise, $\xi$ is overtwisted.}
\end{table}

Let $f$ be the function, specified by this table, which takes a fractional Dehn twist coefficient $c$ and outputs $f(c)=gr(c^+(\xi))$. The theorem below then follows from Proposition \ref{prop:grading}, the grading shift formula in Equation \ref{eqn:grshift}, and the expression for $c_1(\spc_0)^2$ in Equation \ref{eqn:grshift2}.

\begin{theorem}
\label{thm:periodic3}
Suppose $(M,\xi)$ is compatible with a genus one open book $(S_{g,r},\phi)$, where $\phi$ is periodic with fractional Dehn twist coefficients $0\leq k_1/m\leq \dots \leq k_r/m$. Let $I$ be the smallest integer such that $k_I \neq 0$. For $i=I-1,\dots,r-1,$ define $$s_i = -k_r(1/k_r + 1/k_{r-1} + \dots + 1/k_{i+1}).$$ Then $$d_3(\xi) =- f(\frac{k_r}{m})+ \frac{3I-r-4}{4}  -\frac{1}{4}\sum_{j=I-1}^{r-2} \frac{(k_r(j-r) - k_{j+1}s_j)^2}{mk_{j+1}s_js_{j+1}}.$$
\end{theorem}

\section{Capping off and pseudo-Anosov open books}
\label{sec:cappa}
\begin{figure}[!htbp]

\labellist 
\hair 2pt 
\tiny\pinlabel $x$ at 49 74
\pinlabel $x$ at 195 73
\pinlabel $x$ at 342 72
\small\pinlabel $p=3$ at 53 5
\pinlabel $p=4$ at 200 5
\pinlabel $p=5$ at 348 5

\endlabellist 

\begin{center}
\includegraphics[width=11.4cm]{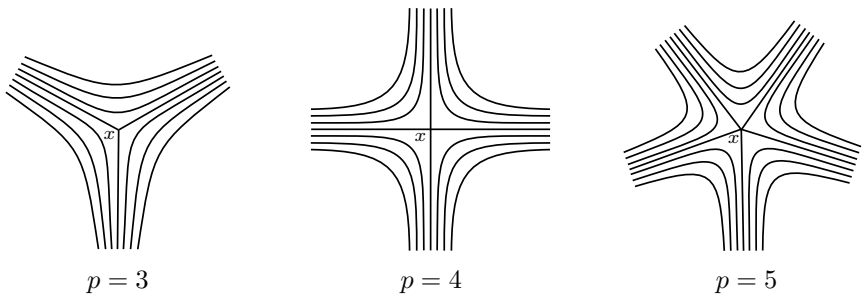}
\caption{\quad Neighborhoods of interior singularities. The singular leaves in each neighborhood are called ``prongs." From left to right, $x$ is a $p$-pronged singularity with $p=3,4,5$. }
\label{fig:sing}
\end{center}
\end{figure}

\begin{figure}[!htbp]

\labellist 
\hair 2pt 
\tiny\pinlabel $x$ at 58 10

\endlabellist 

\begin{center}
\includegraphics[width=4.4cm]{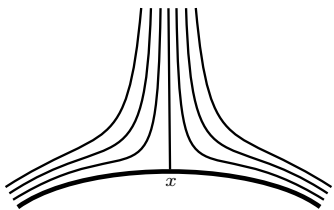}
\caption{\quad A neighborhood of a boundary singularity $x$. The thickened segment represents a portion of the boundary.}
\label{fig:sing2}
\end{center}
\end{figure}

If $\phi$ is a boundary-fixing diffeomorphism of $S$ which is neither periodic nor reducible, then $\phi$ is said to be pseudo-Anosov. In this case (and only in this case), $\phi$ is freely isotopic to a homeomorphism $\phi_0$ for which there exists a transverse pair of singular measured foliations, $(\mathcal{F}_s,\mu_s)$ and $(\mathcal{F}_u,\mu_u)$, of $S$ such that $\phi_0(\mathcal{F}_s,\mu_s)=(\mathcal{F}_s,\lambda\mu_s)$ and $\phi_0(\mathcal{F}_u,\mu_u)=(\mathcal{F}_u,\lambda^{-1}\mu_u)$ for some $\lambda>1$ \cite{th2}. The singularities of $\mathcal{F}_s$ and $\mathcal{F}_u$ which lie in the interior of $S$ are required to be ``$p$-pronged saddles" with $p\geq 3$, as shown in Figure \ref{fig:sing}. Each foliation must have at least one singularity on every boundary component, and each boundary singularity must have a neighborhood of the form shown in Figure \ref{fig:sing2}.

The \emph{fractional Dehn twist coefficient} of $\phi$ around a boundary component of $S$ measures the amount of twisting around this component that takes place in the free isotopy from $\phi$ to $\phi_0$. More precisely, let $x_1,\dots, x_n$ be the singularities of $\mathcal{F}_s$ which lie on some boundary component $B$, labeled in order as one traverses $B$ in the direction specified by its orientation. The map $\phi_0$ permutes these singularities; in fact, we may assume that there exists an integer $k$ for which $\phi_0$ sends $x_i$ to $x_{i+k}$ for all $i$ (where the subscripts are taken modulo $n$). If $H:S\times [0,1]\rightarrow S$ is the free isotopy from $\phi$ to $\phi_0$, and $\beta:B\times [0,1]\rightarrow B\times [0,1]$ is the map which sends $(x,t)$ to $(H(x,t),t)$, then $\beta(x_i \times [0,1])$ is an arc from $(x_i,0)$ to $(x_{i+k},1)$. The fractional Dehn twist coefficient of $\phi$ around $B$ is defined to be the fraction $c \in \mathbb{Q}$, where $c\equiv k/n \text{ modulo }1$ is the number of times that $\beta(x_i\times [0,1])$ wraps around $B\times [0,1]$ (see \cite{hkm1} for more details).

Suppose that $(S_{g,r},\phi)$ is an open book with pseudo-Anosov $\phi$, and let $B_1,\dots,B_r$ denote the boundary components of $S_{g,r}$. Let $c_i$ be the fractional Dehn twists coefficient of $\phi$ around $B_i$, and suppose that $\mathcal{F}_s$ and $\mathcal{F}_u$ are the singular foliations associated to $\phi$. If each of these foliations has $p>1$ singularities on $B_r$ ($\mathcal{F}_s$ and $\mathcal{F}_u$ will have the same number), then they can be extended to transverse singular measured foliations, $\mathcal{F}_s'$ and $\mathcal{F}_u'$, of the surface $S_{g,r-1}$ obtained by capping off the boundary component $B_r$. To see this, remove the leaf corresponding to $B_r$ in each of these foliations, and extend them across the capping disk, creating a $p$-pronged singularity at the center of the disk (if $p=2$, then the foliations extend without singularity over the disk).  The induced diffeomorphism $\phi'$ of the capped off surface $S_{g,r-1}$ is then pseudo-Anosov with associated foliations $\mathcal{F}_s'$ and $\mathcal{F}_u'$. Moreover, $c_i$ is the fractional Dehn twist coefficient of $\phi'$ around the boundary component $B_i$ for $i=1,\dots,r-1$ since this modification took place locally. The requirement that $\mathcal{F}_s$ and $\mathcal{F}_u$ have $p>1$ singularities on $B_r$ is critical in order for this to work; otherwise, there is no obvious way of extending these foliations across the capping disk so that the new interior singularities have $p \geq 3$ prongs. In fact, there \emph{are} pseudo-Anosov diffeomorphisms of $S_{g,r}$ for which the induced diffeomorphism on the capped off surface $S_{g,r-1}$ is \emph{not} pseudo-Anosov. 





In \cite{bald1}, we show that if $\phi$ is a pseudo-Anosov diffeomorphism of $S_{1,1}$ with fractional Dehn twist coefficient less than 1, then $c^+(S_{1,1},\phi)$ is in the image of $U^d$ for all $d \in \mathbb{N}$. The corollary below follows immediately from this fact.

\begin{corollary}
\label{cor:pag1}
Suppose that $(M,\xi)$ is supported by a genus one open book $(S_{1,r},\phi)$ with pseudo-Anosov $\phi$ such that the associated foliations have exactly two singularities on every boundary component of $S_{1,r}$. If any of the fractional Dehn twists coefficients of $\phi$ are less than 1, then $c^+(\xi)$ is in the image of $U^d$ for all $d \in \mathbb{N}$.
\end{corollary}

\begin{remark}
The assumption in Corollary \ref{cor:pag1} that there are \emph{exactly} two singularities on every boundary component is equivalent to the condition that there are \emph{at least} two singularities on each boundary component (and no interior singularities), and is also equivalent to the condition that the foliations associated to $\phi$ are orientable. Finally, note that any open book $(S_{1,r},\phi)$ of the sort considered in the corollary above arises from an Anosov map $\phi_0$ of $S_{1,0}$ by puncturing the torus (creating boundary components) at $r$ fixed points of $\phi_0$.

\end{remark}

\begin{proof}[Proof of Corollary \ref{cor:pag1}]
Suppose the fractional Dehn twist coefficient of $\phi$ around some boundary component is less than 1. After capping off every other boundary component, we obtain an open book $(S_{1,1},\phi')$ where $\phi'$ is pseudo-Anosov with fractional Dehn twist coefficient less than 1 (we may do this since the foliations associated to $\phi$ have more than one singularity on every boundary component). Then, $c^+(S_{1,1},\phi')$ is in the image of $U^d$ for all $d \in \mathbb{N}$. Combined with Theorem \ref{thm:nat}, this proves the corollary.
\end{proof}

In \cite[Theorem 1.1]{hkm2}, Honda, Kazez and Mati{\'c} show (using the taut foliations constructed by Roberts in \cite{roberts1} along with a result of Eliashberg and Thurston \cite{ethur}) that if a pseudo-Anosov diffeomorphism $\phi$ of $S_{1,1}$ has fractional Dehn twist coefficient at least 1, then the contact structure compatible with the open book $(S_{1,1},\phi)$ is weakly symplectically fillable by a filling $W$ with $b_2^+(W)>0$. This prompts the following question.

\begin{question}
\label{ques:symp}
Suppose that $\phi$ is a pseudo-Anosov diffeomorphism of $S_{1,r}$ whose fractional Dehn twist coefficients are all at least 1. Is the contact structure compatible with the open book $(S_{1,r},\phi)$ necessarily weakly symplectically fillable by a filling with $b_2^+(W)>0$?
\end{question} 

In \cite{osz2}, Ozsv{\'a}th and Szab{\'o} show that if $(M,\xi)$ is weakly symplectically fillable by such a filling, and $b_1(M)=0$ (in which case this weak filling may be perturbed to a strong filling \cite{oono}), then there exists some $d \in \mathbb{N}$ for which $c^+(\xi)$ is \emph{not} in the image of $U^d$. The conjecture below follows from this fact, together with Corollary \ref{cor:pag1} and a positive answer to Question \ref{ques:symp}.

\begin{conjecture}
\label{conj:symp}
Suppose that $(M,\xi)$ is supported by a genus one open book $(S_{1,r},\phi)$ with pseudo-Anosov $\phi$ such that the associated foliations have two singularities on every boundary component of $S_{1,r}$. If $b_1(M)=0$, then $\xi$ is strongly symplectically fillable by a filling with $b_2^+(W)>0$ if and only if there exists some $d \in \mathbb{N}$ for which $c^+(\xi)$ is \emph{not} in the image of $U^d$.
\end{conjecture}

One may view Conjecture \ref{conj:symp} as a potential obstruction, via Heegaard Floer homology, to a contact structure being supported by a certain (fairly abundant) type of genus one open book. Even if true, however, this obstruction appears rather cumbersome. One wonders whether there is a more geometric interpretation of the condition that the foliations associated to $\phi$ have two singularities on every boundary component.  Such a condition, combined with a result of the sort proposed in Conjecture \ref{conj:symp} could be helpful in formulating a usable obstruction to $sg(\xi)=1$.

\bibliographystyle{hplain.bst}
\bibliography{References}

\end{document}